\documentclass{elsart}

\usepackage{amssymb}
\usepackage{color}
\usepackage{graphicx}

\usepackage{amsmath}
\numberwithin{equation}{section}

\begin{document}

\begin{frontmatter}
\title{An adaptive $C^{0}$IPG method for the Helmholtz transmission
eigenvalue problem\thanksref{1}}
\thanks[1]{Project supported by the National Natural Science Foundation of China (Grant No.11561014).}
\author{Hao Li},\ead{lihao365@126.com}
\author{ Yidu Yang\corauthref{2}}\ead{ydyang@gznu.edu.cn}
\corauth[2]{Corresponding author}
\address{School of Mathematical Sciences, Guizhou Normal University,\\ Guiyang 550001, China}

\begin{abstract}
The interior penalty methods using $C^{0}$ Lagrange
elements ($C^{0}$IPG) developed in the last decade  for the fourth
order problems are an interesting topic in academia at present. In this
paper, we discuss the adaptive
fashion of $C^{0}$IPG method for the Helmholtz transmission eigenvalue problem.
We give the a posteriori error indicators for primal and dual eigenfunctions, and prove their reliability and efficiency.
We also give the a posteriori error indicator for eigenvalues and design a $C^{0}$IPG adaptive algorithm.
Numerical experiments show that this algorithm is efficient and can get the optimal convergence rate.
\end{abstract}
\end{frontmatter}
\section{Introduction}
\label{intro}
\indent The transmission eigenvalues can be used to obtain estimates for
the material properties of the scattering object
\cite{cakoni1,cakoni2,sun0}, and
have theoretical importance in the uniqueness and reconstruction in
inverse scattering theory \cite{colton1}. In recent years, the computation of transmission
eigenvalues has attracted the attention of many researchers. The
first numerical treatment of the transmission eigenvalue problem
appeared in \cite{colton2} where three finite element methods,
including the Argyris, continuous and mixed methods,
are proposed for the Helmholtz transmission eigenvalues, and has
been further developed  by
\cite{an1,cakoni5,geng,han,ji,ji2,monk,sun1,su4,yang3,yang2,yang4,han2,Kleefeld} et al..\\
\indent $C^{0}$ interior penalty Galerkin ($C^{0}$IPG)
method, developed in the last decade \cite{engel,brenner1}, is a new class of Galerkin
methods for fourth order problems.
The researches for $C^{0}$IPG methods have been an interesting topic in academia at present.
There exist many researches for fourth order elliptic equations (see \cite{engel,brenner1,brenner3,gudi,ji3}) and
  for eigenvalue problems
(see \cite{geng,yang5,wells,brenner2,brennerADD,hao}) by
$C^{0}$IPG methods.\\
The a posteriori error estimates and adaptive finite element methods are always the main
streams of scientific and engineering computing. The idea of the a posteriori error estimates
was first introduced by Babuska and Rheinboldt \cite{babuska1} in 1978. Up to now, many excellent works have been summarized in the books
such as  \cite{ainsworth,verfurth,shi}. And a posteriori error estimates of residual type of $C^{0}$IPG method of fourth order elliptic equations also
have been summarized in \cite{brenner1}.\\
\indent Inspired by the works mentioned above,
in this paper,
based on the weak formulation proposed in \cite{yang2,yang4},
we propose a new $C^{0}$IPG discrete scheme (see (\ref{s2.25}))
and discuss
the a posteriori error estimates and adaptive
algorithm of $C^{0}$IPG method for the Helmholtz transmission eigenvalue
problem.
We give the a posteriori error indicators for primal and dual eigenfunctions and
eigenvalues. We prove that the indicators for both primal and
dual eigenfunctions are reliable and efficient, and analyze the reliability of the indicator for eigenvalues.
Based on the given indicators, we design an adaptive algorithm.
Numerical experiments show that this algorithm is efficient and can get the optimal convergence rate.
Compared with adaptive $C^1$ conforming finite element algorithm in \cite{han},
the adaptive $C^{0}$IPG algorithm is simpler to be constructed and implemented numerically.\\
 \indent In
this paper, regarding the basic theory of  finite
element methods, we refer to \cite{shi,babuska,brenner,ciarlet,oden}.\\
\indent Throughout this paper, the letter $C$ (with or without
subscripts) denotes a positive constant independent of mesh size
$h$, which may not be the same constant in different places. For
simplicity, we use the symbol $a\lesssim b$ to mean that $a\leq C
b$ and the symbol $a\approx b$ to mean $a\lesssim b\lesssim a$.

\section{A $C^{0}$IPG discrete scheme}

 \indent Consider the Helmholtz transmission eigenvalue
problem: Find $k\in \mathbb{C}$, $w, \sigma\in L^{2}(\Omega)$,
$w-\sigma\in H^{2}(\Omega)$ such that
\begin{eqnarray}\label{s2.1}
&&\Delta w+k^{2}nw=0,~~~in~\Omega,\\\label{s2.2}
 &&\Delta
\sigma+k^{2}\sigma=0,~~~in~ \Omega,\\\label{s2.3}
 &&w-\sigma=0,~~~on~ \partial
\Omega,\\\label{s2.4}
 &&\frac{\partial w}{\partial \gamma}-\frac{\partial
\sigma}{\partial \gamma}=0,~~~ on~\partial \Omega,
\end{eqnarray}
where $\Omega\subset \mathbb{R}^{d}~ (d=2,3)$ is a bounded simply
connected inhomogeneous medium,
 $\gamma$ is the unit outward normal to $\partial \Omega$ and the index of refraction $n=n(x)$
is positive.\\
\indent Let $W^{s,p}(\Omega)$ denote the usual Sobolev space with
norm $\|\cdot\|_{s,p}$,
 $H^{s}(\Omega)=W^{s,2}(\Omega)$, and $\|\cdot\|_{s,2}=\|\cdot\|_{s}$,
$H^{0}(\Omega)=L^{2}(\Omega)$ with the inner product
$(u,v)_{0}=\int_{\Omega}u\overline{v}dx$. Denote $
H_{0}^{2}(\Omega)=\{v\in H^{2}(\Omega): v|_{\partial
\Omega}=\frac{\partial v}{\partial \gamma}|_{\partial \Omega}=0\}. $
 Let $H^{-1}(\Omega)$ be the
``negative space" with  norm
$\|v\|_{-1}$.\\
 Define Hilbert
space $\mathbf{H}=H_{0}^{2}(\Omega)\times L^{2}(\Omega)$ with norm
$\|(v,z)\|_{\mathbf{H}}=\|v\|_{2}+\|z\|_{0}$, and define
$\mathbf{H}^{1}=H_0^{1}(\Omega)\times H^{-1}(\Omega)$ with norm
$\|(v,z)\|_{\mathbf{H}^{1}}=\|v\|_{1}+\|z\|_{-1}$.\\
Since $L^{2}(\Omega)\hookrightarrow H^{-1}(\Omega)$ compactly and
$H^{2}(\Omega)\hookrightarrow
H^{1}(\Omega)$ compactly, $\mathbf{H}\hookrightarrow\mathbf{H}^{1}$ compactly. \\
\indent In this paper, we suppose that $n\in W^{1,\infty}(\Omega)$
satisfying the following condition
\begin{eqnarray*}
1+\delta\leq n(x)~in~ \Omega,
\end{eqnarray*}
for some constant $\delta>0$. And the argument is the same if $0<
n(x)\leq 1-\varrho$~in~$\Omega~(\varrho>0)$
holds.\\

\indent From \cite{cakoni3,rynne} we know that the problem
 (\ref{s2.1})-(\ref{s2.4}) can be written as the following equivalent
weak formulation: Find $k\in \mathbb{C}$, $u\in
H_{0}^{2}(\Omega)$ such that
\begin{eqnarray*}
&&(\frac{1}{n-1}\Delta u,\Delta v )_{0}=k^{2}(\nabla u,
\nabla(\frac{n}{n-1}v)_{0}+k^{2}(\nabla (\frac{1}{n-1}u),
\nabla v)_{0}\nonumber\\
&&~~~~~~-k^{4}(\frac{n}{n-1}u, v)_{0},~~~\forall v
\in~H_{0}^{2}(\Omega).
\end{eqnarray*}
Introduce an auxiliary variable $ \omega=k^{2}u$, and let
$\lambda=k^{2}$, then we arrive at a linear weak formulation (see
\cite{yang2,yang4}): Find $\lambda\in \mathbb{C}$, $(u,\omega)\in
\mathbf{H}\setminus \{0\}$ such that
\begin{eqnarray}\label{s2.5}
A((u,\omega),(v,z)) =\lambda B((u,\omega),(v,z)),~~~\forall (v,z)\in
\mathbf{H},
\end{eqnarray}
where
\begin{eqnarray}\label{s2.6}
A((u,\omega),(v,z))=((\frac{1}{n-1}-\mu)\Delta u, \Delta
v)_{0}+\mu\int\limits_{\Omega}D^{2}u:D^{2}\bar{v}dx+(\omega, z)_{0}
\end{eqnarray}
with constant $\mu>0$, $\frac{1}{n-1}-\mu\geq 0$, and
\begin{eqnarray*}
B((u,\omega),(v,z))=(\nabla(\frac{1}{n-1}u), \nabla v)_{0}+(\nabla
u, \nabla(\frac{n}{n-1}v))_{0}-(\omega, \frac{n}{n-1}v
)_{0}+(u,z)_{0}.
\end{eqnarray*}
\indent  It is obvious that $A(\cdot,\cdot)$ is a selfadjoint,
continuous sesquilinear form on $\mathbf{H}\times \mathbf{H}$,
\begin{eqnarray}\label{s2.7}
A((v,z),(v,z))\gtrsim \|(v,z)\|_{\mathbf{H}}^{2},
\end{eqnarray}
and for any given $(f,g)\in
\mathbf{H}^{1}$, $B((f,g), (v, z))$ is a continuous linear form on
$\mathbf{H}$,
\begin{eqnarray}\label{s2.8}
|B((f,g), (v, z))|\lesssim
\|(f,g)\|_{\mathbf{H}^{1}}\|(v,z)\|_{\mathbf{H}^{1}},~~~\forall (v, z)\in \mathbf{H}^{1}.
\end{eqnarray}
 \indent We use $A(\cdot,\cdot)$ and
$\|\cdot\|_{A}=A(\cdot,\cdot)^{\frac{1}{2}}$ as an inner product and
norm on $\mathbf{H}$, respectively.\\
\indent The source problem associated with (\ref{s2.5}) is as
follows: Find $(\psi,\varphi)\in \mathbf{H}$ such that
\begin{eqnarray}\label{s2.9}
A((\psi,\varphi),(v,z)) =B((f,g),(v,z)),~~~\forall (v,z)\in
\mathbf{H}.
\end{eqnarray}
From Lax-Milgram theorem we know that (\ref{s2.9}) has one and only
one solution. Therefore, we define the corresponding solution
operator $T: \mathbf{H}^{1}\to \mathbf{H}$ by
\begin{eqnarray}\label{s2.10}
A(T(f,g),(v,z)) =B((f,g), (v,z)),~~~\forall (v,z)\in \mathbf{H}.
\end{eqnarray}
Then (\ref{s2.5}) has the equivalent operator form:
\begin{eqnarray}\label{s2.11}
T(u,\omega)=\frac{1}{\lambda}(u,\omega).
\end{eqnarray}

\indent From (\ref{s2.10}) we have
\begin{eqnarray}\label{s2.12}
\|T(f,g)\|_{\mathbf{H}}\lesssim \|(f,g)\|_{\mathbf{H}^{1}},~~~\forall (f,g) \in \mathbf{H}^{1}.
\end{eqnarray}
Thus we know that $T: \mathbf{H}\to \mathbf{H}$ is compact,
and $T: \mathbf{H}^{1}\to
\mathbf{H}^{1}$ is compact.\\
\indent Consider the dual problem of (\ref{s2.5}): Find
$\lambda^{*}\in \mathbb{C}$, $(u^{*},\omega^{*})\in
\mathbf{H}\setminus \{0\}$ such that
\begin{eqnarray}\label{s2.13}
A((v,z), (u^{*},\omega^{*})) =\overline{\lambda^{*}} B((v,z),
(u^{*},\omega^{*})),~~~\forall (v,z)\in \mathbf{H}.
\end{eqnarray}
\indent The source problem associated with (\ref{s2.13}) is as
follows: Find $(\psi^{*},\varphi^{*})\in \mathbf{H}$ such that
\begin{eqnarray}\label{s2.14}
A((v,z), (\psi^{*},\varphi^{*})) =B((v,z), (f,g)),~~~\forall
(v,z)\in \mathbf{H}.
\end{eqnarray}
Define the corresponding solution operator $T^{*}:
\mathbf{H}^{1}\to \mathbf{H}$ by
\begin{eqnarray}\label{s2.15}
A((v,z), T^{*}(f,g)) =B((v,z), (f,g)),~~~\forall (v,z)\in
\mathbf{H}.
\end{eqnarray}
Then (\ref{s2.13}) has the equivalent operator form:
\begin{eqnarray}\label{s2.16}
T^{*}(u^{*},\omega^{*})=\lambda^{*-1}(u^{*},\omega^{*}).
\end{eqnarray}
\indent From (\ref{s2.10}) and (\ref{s2.15}) we know that
$T^{*}$ is the adjoint operator of $T$ in the sense of inner product
$A(\cdot,\cdot)$. So the primal and dual eigenvalues are connected
via
$\lambda=\overline{\lambda^{*}}$ (see \cite{yang2}).\\

\indent Denote
$$\mathbb{S}=(\frac{2d}{1+d}, 2].$$

\indent We need the following regularity assumption:\\
\indent {\bf $R(\Omega)$.}~~{\em For any $\xi\in H^{-1}(\Omega)$,
there exists $\psi\in W^{3,p_{0}}(\Omega)$ satisfying
\begin{eqnarray*}
\Delta(\frac{1}{n-1}\Delta
\psi)=\xi,~~~in~\Omega,~~~\psi=\frac{\partial \psi}{\partial
\gamma}=0~~~ on~\partial \Omega,
\end{eqnarray*}
and
\begin{eqnarray}\label{s2.17}
\|\psi\|_{3,p_{0}}\leq C_{\Omega} \|\xi\|_{-1},
\end{eqnarray}
where $p_{0}\in \mathbb{S}$, $C_{\Omega}$ denotes the  prior constant
dependent on the $n(x)$ and $\Omega$ but independent of the
right-hand
side $\xi$ of the equation.}\\
\indent Let $\pi_{h}$ be a shape-regular mesh, for any element
$\kappa\in\pi_{h}$, let $h_{\kappa}$ denote diameter of $\kappa$,
$h=\max_{\kappa\in\pi_{h}}h_{\kappa}$. And let
$$S^{h}=\{v\in C(\bar{\Omega})\cap H_{0}^{1}(\Omega):
v|_{\kappa}\in P_{m}, \forall \kappa\in\pi_{h}\},$$ where
$P_{m}$ is the set of all polynomials in $d$
variables of degree $\leq m (m\geq2)$.
Let
$\mathbf{H}_{h}=S^{h}\times S^{h}$. Then $\mathbf{H}_{h}\subset
\mathbf{H}^{1}$ but $\mathbf{H}_{h}\not\subset \mathbf{H}$.\\
\indent Let $p\in \mathbb{S}$, from the trace theorem with scaling
we have the following trace inequality:
\begin{eqnarray}\label{s2.18}
&&\int\limits_{\partial\kappa} |w|^{2}ds\lesssim
h_{\kappa}^{d-\frac{2
d}{p}-1}\|w\|_{0,p,\kappa}^{2}+h_{\kappa}^{2+d-\frac{2
d}{p}-1}|w|_{1,p,\kappa}^{2},~~~\forall \kappa\in \pi_{h}.
\end{eqnarray}
\indent Let $\mathcal {E}$ denote the set of all $(d-1)$-faces in
$\pi_{h}$ $(d=2,3)$.
We decompose $\mathcal {E}=\mathcal {E}^{i}\cup\mathcal
{E}^{b}$ where $\mathcal {E}^{i}$ and $\mathcal
{E}^{b}$ refer to interior faces and faces on the
boundary $\partial \Omega$, respectively. For each $\ell\in \mathcal
{E}^{i}$, we choose an arbitrary unit normal vector
$\gamma_{\ell}$ and denote the two triangles sharing this face by
$\kappa_{-}$ and $\kappa_{+}$, where $\gamma_{\ell}$ points from
$\kappa_{-}$ to $\kappa_{+}$. We set the jump and average on $\ell$
as
\begin{eqnarray}\label{s2.19}
&&[[\frac{\partial v}{\partial \gamma_{\ell}}]]=\bigtriangledown
(v|_{\kappa_{+}})\cdot \gamma_{\ell} -\bigtriangledown
(v|_{\kappa_{-}})\cdot \gamma_{\ell},\\\label{s2.20}
 &&\{\{(\frac{1}{n-1}-\mu)\Delta
v\}\}=\frac{1}{2}((\frac{1}{n-1}-\mu)\Delta
v|_{\kappa_{-}}+(\frac{1}{n-1}-\mu)\Delta
v|_{\kappa_{+}}),\\\label{s2.21}
&&\{\{\frac{\partial^{2}v}{\partial\gamma_{\ell}^{2}}\}\}=\frac{1}{2}(\frac{\partial^{2}v}{\partial\gamma_{\ell}^{2}}|_{\kappa_{-}}
+\frac{\partial^{2}v}{\partial\gamma_{\ell}^{2}}|_{\kappa_{+}})
\end{eqnarray}
with
$\frac{\partial^{2}v}{\partial\gamma_{\ell}^{2}}=\gamma_{\ell}\cdot(D^{2}v)\gamma_{\ell}$.\\
For any $\ell\in \mathcal {E}^{b}$ which is a face of
$\kappa$, we take $\gamma_{\ell}$ to be the unit normal vector
pointing towards the outside of $\Omega$ and set
\begin{eqnarray}\label{s2.22}
&&[[\frac{\partial v}{\partial \gamma_{\ell}}]]=-\gamma_{\ell}\cdot
\bigtriangledown (v|_{\kappa}),\\\label{s2.23}
&&\{\{(\frac{1}{n-1}-\mu)\Delta v\}\}=(\frac{1}{n-1}-\mu)\Delta
v|_{\kappa},~~~\{\{\frac{\partial^{2}v}{\partial\gamma_{\ell}^{2}}\}\}=\frac{\partial^{2}v}{\partial\gamma_{\ell}^{2}}|_{\kappa}.
\end{eqnarray}
Define piecewise Sobolev space
$$W^{3,p}(\Omega, \pi_{h})=\{v\in C(\bar{\Omega})\cap H_{0}^{1}(\Omega): v|_{\kappa}\in W^{3,p}(\kappa)~\forall \kappa\in\pi_{h}\}~~~(p\in\mathbb{S}). $$
Referring \cite{geng,brenner1,yang5}, we define
\begin{eqnarray}\label{s2.24}
&&A_{h}((u,\omega),(v,z))
=\sum\limits_{\kappa\in\pi_{h}}\int\limits_{\kappa}
(\frac{1}{n-1}-\mu)\Delta u \Delta \overline{v}dx+
\mu\int\limits_{\kappa} D^{2}u: D^{2} \overline{v}dx
\nonumber\\
&&~~~~~~+ \sum\limits_{\ell\in \mathcal {E}}\int\limits_{\ell}
\{\{(\frac{1}{n-1}-\mu)\Delta u\}\}[[\frac{\partial
\overline{v}}{\partial
\gamma_{\ell}}]]+\{\{(\frac{1}{n-1}-\mu)\Delta v\}\}[[\frac{\partial
\overline{u}}{\partial \gamma_{\ell}}]]ds\nonumber\\
&&~~~~~~+ \mu\sum\limits_{\ell\in \mathcal {E}}\int\limits_{\ell}
\{\{\frac{\partial^{2}u}{\partial\gamma_{\ell}^{2}}\}\}[[\frac{\partial
\overline{v}}{\partial
\gamma_{\ell}}]]+\{\{\frac{\partial^{2}v}{\partial\gamma_{\ell}^{2}}\}\}[[\frac{\partial
\overline{u}}{\partial \gamma_{\ell}}]]ds\nonumber\\
&&~~~~~~ + \sigma\sum\limits_{\ell\in \mathcal
{E}}\frac{1}{\hat{\ell}}\int\limits_{\ell}[[\frac{\partial
u}{\partial \gamma_{\ell}}]] [[\frac{\partial \overline{v}}{\partial
\gamma_{\ell}}]]ds+
\sum\limits_{\kappa\in\pi_{h}}\int\limits_{\kappa}\omega
\overline{z}dx,
\end{eqnarray}
where $\sigma>1$ is the penalty parameter,
and $\hat{\ell}=h_{\ell}$ is the diameter of $\ell$.\\
\indent We give the following $C^{0}$IPG discrete scheme of
(\ref{s2.5}): Find $\lambda_{h}\in \mathbb{C}$,
$(u_{h},\omega_{h})\in \mathbf{H}_{h}\setminus \{0\}$ such that
\begin{eqnarray}\label{s2.25}
A_{h}((u_{h},\omega_{h}),(v,z)) =\lambda_{h}
B((u_{h},\omega_{h}),(v,z)),~~~\forall (v,z)\in \mathbf{H}_{h}.
\end{eqnarray}
\indent We define the mesh-dependent norms $\|\cdot\|_{h}$ and
$\||\cdot\||_{h}$ on $W^{3,p}(\Omega, \pi_{h})\times L^{2}(\Omega)$
as
\begin{eqnarray}\label{s2.26}
&&\|(u,\omega)\|_{h}^{2} =\sum\limits_{\kappa\in\pi_{h}} \|
u\|_{2,\kappa}^{2}
 + \sigma\sum\limits_{\ell\in \mathcal
{E}}\frac{1}{\hat{\ell}} \|[[\frac{\partial u}{\partial
\gamma_{\ell}}]]\|_{0,\ell}^{2} +
\sum\limits_{\kappa\in\pi_{h}}\|\omega\|_{0,\kappa}^{2},\\\label{s2.27}
&&\||(u,\omega)\||_{h}^{2}=\|(u,\omega)\|_{h}^{2}
+\frac{1}{\sigma}\sum\limits_{\ell\in \mathcal {E}} \|\{\{\Delta
u\}\}\|_{0,\ell}^{2}\hat{\ell}+ \frac{1}{\sigma}\sum\limits_{\ell\in
\mathcal {E}} \|\{\{\frac{\partial^{2} u}{\partial
\gamma_{\ell}^{2}} \}\}\|_{0,\ell}^{2}\hat{\ell}.
\end{eqnarray}
By the trace inequality (\ref{s2.18}) with $p=2$ and  the inverse
estimates we have
\begin{eqnarray}\label{s2.28}
\|\Delta v\|_{0,\ell}\lesssim \hat{\ell}^{-\frac{1}{2}}\|
v\|_{2,\kappa},~~ \|\{\{\frac{\partial^{2} u}{\partial
\gamma_{\ell}^{2}} u\}\}\|_{0,\ell}\lesssim
\hat{\ell}^{-\frac{1}{2}}\| v\|_{2,\kappa}, ~~~\forall v\in S^{h}.
\end{eqnarray}
So on $\mathbf{H}_{h}$ the two norms $\|\cdot\|_{h}$ and
$\||\cdot\||_{h}$ are equivalent.\\
For any $(u,\omega),(v,z)\in W^{3,p}(\Omega, \pi_{h})\times
L^{2}(\Omega)$, by the Schwartz inequality we can deduce
\begin{eqnarray}\label{s2.29}
&&|A_{h}((u,\omega),(v,z))| \lesssim \sum\limits_{\kappa\in\pi_{h}}
\|\Delta u\|_{0,\kappa}\| \Delta \overline{v}\|_{0,\kappa}+
\sum\limits_{\kappa\in\pi_{h}} |u|_{2,\kappa}|
\overline{v}|_{2,\kappa}
\nonumber\\
&&+ \sum\limits_{\ell\in \mathcal {E}}(
\sqrt{\frac{\hat{\ell}}{\sigma}}\|\{\{\Delta
u\}\}\|_{0,\ell}\sqrt{\frac{\sigma}{\hat{\ell}}}\|[[\frac{\partial
\overline{v}}{\partial
\gamma_{\ell}}]]\|_{0,\ell}+\sqrt{\frac{\hat{\ell}}{\sigma}}\|\{\{\Delta
v\}\}\|_{0,\ell}\sqrt{\frac{\sigma}{\hat{\ell}}}\|[[\frac{\partial
\overline{u}}{\partial \gamma_{\ell}}]]\|_{0,\ell})\nonumber\\
&&+ \sum\limits_{\ell\in \mathcal {E}}(
\sqrt{\frac{\hat{\ell}}{\sigma}}\|\{\{\frac{\partial^{2} u}{\partial
\gamma_{\ell}^{2}}
\}\}\|_{0,\ell}\sqrt{\frac{\sigma}{\hat{\ell}}}\|[[\frac{\partial
\overline{v}}{\partial
\gamma_{\ell}}]]\|_{0,\ell}+\sqrt{\frac{\hat{\ell}}{\sigma}}\|\{\{\frac{\partial^{2}
v}{\partial \gamma_{\ell}^{2}}
\}\}\|_{0,\ell}\sqrt{\frac{\sigma}{\hat{\ell}}}\|[[\frac{\partial
\overline{u}}{\partial \gamma_{\ell}}]]\|_{0,\ell})\nonumber\\
&&+ \sigma\sum\limits_{\ell\in \mathcal
{E}}\frac{1}{\sqrt{\hat{\ell}}}\|[[\frac{\partial u}{\partial
\gamma_{\ell}}]]\|_{0,\ell}\frac{1}{\sqrt{\hat{\ell}}}
\|[[\frac{\partial \overline{v}}{\partial
\gamma_{\ell}}]]\|_{0,\ell}+
\sum\limits_{\kappa\in\pi_{h}}\|\omega\|_{0,\kappa}
\|\overline{z}\|_{0,\kappa}\nonumber\\
 &&\lesssim
\||(u,\omega)\||_{h}\||(v,z)\||_{h}.
\end{eqnarray}
And for any $(u_{h},\omega_{h}), (v,z)\in \mathbf{H}_{h}$, we have
\begin{eqnarray}\label{s2.30}
|A_{h}((u_{h},\omega_{h}),(v,z))|\leq
C\|(u_{h},\omega_{h})\|_{h}\|(v,z)\|_{h}.
\end{eqnarray}
And referring \cite{geng,gudi}, when $\sigma$ is large enough, by
(\ref{s2.28}) and the Young inequality we deduce
\begin{eqnarray}\label{s2.31}
&&A_{h}((u_{h},\omega_{h}),(u_{h},\omega_{h})) \geq C_{1}
\sum\limits_{\kappa\in\pi_{h}}( \|\Delta
u_{h}\|_{0,\kappa}^{2}+\|u_{h}\|_{2,\kappa}^{2})\nonumber\\
&&~~~~~~ -\sqrt{C_{1}} (\sum\limits_{\kappa\in\pi_{h}} \|\Delta
u_{h}\|_{0,\kappa}^{2})^{\frac{1}{2}}
\frac{C}{\sqrt{C_{1}}}(\sum\limits_{\ell\in \mathcal
{E}}\frac{1}{\hat{\ell}} \|[[\frac{\partial u_{h}}{\partial
\gamma_{\ell}}]]\|_{0,\ell}^{2})^{\frac{1}{2}}\nonumber\\
&&~~~~~~ -\sqrt{C_{1}} (\sum\limits_{\kappa\in\pi_{h}} \|
u_{h}\|_{2,\kappa}^{2})^{\frac{1}{2}}
\frac{C}{\sqrt{C_{1}}}(\sum\limits_{\ell\in \mathcal
{E}}\frac{1}{\hat{\ell}} \|[[\frac{\partial u_{h}}{\partial
\gamma_{\ell}}]]\|_{0,\ell}^{2})^{\frac{1}{2}}\nonumber\\
&&~~~~~~ +\sigma\sum\limits_{\ell\in \mathcal
{E}}\frac{1}{\hat{\ell}} \|[[\frac{\partial u_{h}}{\partial
\gamma_{\ell}}]]\|_{0,\ell}^{2}
+\sum\limits_{\kappa\in\pi_{h}}\|\omega_{h}\|_{0,\kappa}^{2}\nonumber\\
&&~~~\geq\frac{C_{1}}{2} \sum\limits_{\kappa\in\pi_{h}} \|\Delta
u_{h}\|_{0,\kappa}^{2}
 +( \frac{\sigma}{2}-\frac{1}{2}\frac{C^{2}}{C_{1}})\sum\limits_{\ell\in \mathcal
{E}}\frac{1}{\hat{\ell}} \|[[\frac{\partial u_{h}}{\partial
\gamma_{\ell}}]]\|_{0,\ell}^{2}\nonumber\\
&&~~~~~~+\frac{C_{1}}{2} \sum\limits_{\kappa\in\pi_{h}} \|
u_{h}\|_{2,\kappa}^{2}
 +( \frac{\sigma}{2}-\frac{1}{2}\frac{C^{2}}{C_{1}})\sum\limits_{\ell\in \mathcal
{E}}\frac{1}{\hat{\ell}} \|[[\frac{\partial u_{h}}{\partial
\gamma_{\ell}}]]\|_{0,\ell}^{2} +
\sum\limits_{\kappa\in\pi_{h}}\|\omega_{h}\|_{0,\kappa}^{2}\nonumber\\
&&~~~\gtrsim \|(u_{h},\omega_{h})\|_{h}^{2},~~~\forall
(u_{h},\omega_{h})\in \mathbf{H}_{h}.
\end{eqnarray}
\indent Consider the $C^{0}$IPG discrete scheme of (\ref{s2.9}): Find
$(\psi_{h},\varphi_{h})\in \mathbf{H}_{h}$ such that
\begin{eqnarray}\label{s2.32}
A_{h}((\psi_{h},\varphi_{h}),(v,z)) =B((f,g), (v,z)),~~~\forall
(v,z)\in \mathbf{H}_{h}.
\end{eqnarray}
\indent We introduce the corresponding solution operator: $T_{h}:
\mathbf{H}^{1}\to \mathbf{H}_{h}$:
\begin{eqnarray}\label{s2.33}
A_{h}(T_{h}(f,g),(v,z)) =B((f,g), (v,z)),~~~\forall (v,z)\in
\mathbf{H}_{h}.
\end{eqnarray}
Then (\ref{s2.25}) has the operator form:
\begin{eqnarray}\label{s2.34}
T_{h}(u_{h},\omega_{h})=\frac{1}{\lambda_{h}}(u_{h},\omega_{h}).
\end{eqnarray}
\indent The $C^{0}$IPG discrete scheme of (\ref{s2.13}) is given by:
 Find
$\lambda_{h}^{*} \in \mathbb{C}$, $(u_{h}^{*},\omega_{h}^{*})\in
\mathbf{H}_{h}\setminus \{0\}$ such that
\begin{eqnarray}\label{s2.35}
A_{h}((v,z), (u_{h}^{*},\omega_{h}^{*})) =\overline{\lambda_{h}^{*}}
B((v,z), (u_{h}^{*},\omega_{h}^{*})),~~~\forall (v,z)\in
\mathbf{H}_{h}.
\end{eqnarray}
Define the solution operator $T_{h}^{*}: \mathbf{H}^{1}\to
\mathbf{H}_{h}$ satisfying
\begin{eqnarray}\label{s2.36}
A_{h}((v,z), T_{h}^{*}(f,g))=B((v,z), (f,g)),~~~\forall~(v,z)\in
\mathbf{H}_{h}.
\end{eqnarray}
Thus (\ref{s2.35}) has the following equivalent operator form:
\begin{eqnarray}\label{s2.37}
T_{h}^{*}(u_{h}^{*},\omega_{h}^{*})=\lambda_{h}^{*-1}(u_{h}^{*},\omega_{h}^{*}).
\end{eqnarray}
\indent It can be proved that $T_{h}^{*}$ is the adjoint operator of
$T_{h}$ in the sense of inner product $A_{h}(\cdot,\cdot)$. In fact,
$\forall (u,\omega),(v, z)\in \mathbf{H}_{h}$, from (\ref{s2.33})
and (\ref{s2.36}) we have
\begin{eqnarray*}
A_{h}(T_{h}(u,\omega),(v, z))=B((u,\omega),(v,
z))=A_{h}((u,\omega),T_{h}^{*}(v, z)).
\end{eqnarray*}
Hence, the primal and dual eigenvalues are connected via
$\lambda_{h}=\overline{\lambda_{h}^{*}}$.\\
\indent In this paper, we suppose that $\{\lambda_{j}\}$ and
$\{\lambda_{j,h}\}$ are    enumerations of the eigenvalues of
(\ref{s2.5}) and (\ref{s2.25}) respectively according to the same
sort rule,  each repeated as many times as its multiplicity, and
$\lambda=\lambda_{i}$ is the $i$th eigenvalue with the algebraic
multiplicity $q$ and the ascent $\alpha$,
$\lambda_{i}=\lambda_{i+1}=\cdots,\lambda_{i+q-1}$, and $\lambda_{h}=\lambda_{i,h}$. When
$\|T_{h}-T\|_{\mathbf{H}^{1}}\to 0$, $q$ eigenvalues
$\lambda_{i,h},\cdots,\lambda_{i+q-1,h}$ of (\ref{s2.25}) will
converge to $\lambda$.\\
\indent Let $E$ be the spectral projection associated with $T$ and
$\lambda$, then $ran(E)=null((\lambda^{-1}-T)^{\alpha})$ is the
space of generalized eigenfunctions associated with $\lambda$ and
$T$, where $ran$ denotes the range and $null$ denotes the null
space. Let $E_{h}
$ be the spectral projection associated with
$T_{h}$ and the eigenvalues
$\lambda_{i,h},\cdots,\lambda_{i+q-1,h}$, then $ ran(E_{h}) $ is the
space spanned by all generalized eigenfunctions corresponding to all
eigenvalues $\lambda_{i,h},\cdots,\lambda_{i+q-1,h}$. In view of the
adjoint problem (\ref{s2.13}) and (\ref{s2.35}), the definitions of
$E^{*}$, $ran(E^{*})$, $E_{h}^{*}$ and $ran(E_{h}^{*})$ are
analogous to $E$,
$ran(E)$, $E_{h}$ and $ran(E_{h})$  (see \cite{babuska}).\\

The error estimate of the C$^0$IPG method for eigenvalue problems is based on the error estimate of the C$^0$IPG method for the corresponding source problems.
Next using argument as in \cite{yang5} we well prove the a priori error estimates for the source problem (\ref{s2.9}).\\
From Lemma 3.1 in  \cite{yang5} we known that
(\ref{s2.9}) admits a
 unique solution
$(\psi,\varphi)\in (W^{3,p_{0}}(\Omega)\cap H_{0}^{2}(\Omega))\times
H_{0}^{1}(\Omega)$ and
\begin{eqnarray}\label{s2.38r}
\|(\psi,\varphi)\|_{W^{3,p_{0}}\times
H_{0}^{1}}\leq C_{R}\|(f,g)\|_{\mathbf{H}^{1}},~~~\forall (f,g)\in
\mathbf{H}^{1},
\end{eqnarray}
where $p_{0}\in \mathbb{S}$, $C_{R}$ denotes the  prior constant.\\
\indent Denote $A((u,\omega), (v,z))\equiv a(u,v)+(\omega, z)_{0}$, $A_{h}((u,\omega), (v,z))\equiv a_{h}(u,v)+(\omega, z)_{0}$,
$\|(u,\omega)\|_{h}\equiv\|u\|_{h}^{2}+\|\omega\|_{0}^{2}$, $\||(u,\omega)\||_{h}\equiv\||u\||_{h}^{2}+\|\omega\|_{0}^{2}$,
$B'(f, v)=\int\limits_{\Omega}\nabla
f\cdot\nabla \bar{v}dx$. Define the
auxiliary operator $K: H_{0}^{1}(\Omega)\to  H_{0}^{2}(\Omega)$ by
\begin{eqnarray}\label{s2.39r}
a(Kf, v)=B'(f, v),~~~\forall v\in H_{0}^{2}(\Omega).
\end{eqnarray}
Then for any $f\in H_{0}^{1}(\Omega)$ it is valid that
$Kf\in W^{3,p_{0}}(\Omega)$ and
\begin{eqnarray}\label{s2.40r}
&&\|Kf\|_{3,p_{0}}\lesssim\|f\|_{1
}.
\end{eqnarray}
Referring (3.7)-(3.9) in  \cite{yang5}  we can deduce
\begin{eqnarray}\label{s2.41r}
&&A_{h}((\psi,\varphi), (v,z))=B((f,g), (v,z)),~~~\forall (v,z)\in
\mathbf{H}_{h},\\\label{s2.42r}
&&a_{h}(Kf, v)=B'(f, v),~~~\forall v\in
S^{h}.
\end{eqnarray}
From (\ref{s2.41r}) and (\ref{s2.32}) we get
\begin{eqnarray}\label{s2.43r}
A_{h}((\psi,\varphi)-(\psi_{h},\varphi_{h}), (v,z))=0,~~~\forall
(v,z)\in \mathbf{H}_{h}.
\end{eqnarray}
 \indent Define the operator
 $$I_{h}(\psi,\varphi)=(I_{h}^{1}\psi, I_{h}^{2}\varphi),$$
where $I_{h}^{1}: H_{0}^{1}(\Omega)\cap C^{0}(\overline{\Omega})\to
S^{h}$ is the Lagrange nodal interpolation operator and $I_{h}^{2}: L^{2}(\Omega)\to
S^{h}$  is defined by
\begin{eqnarray*}
(\varphi-I_{h}^{2}\varphi, z)_{0}=0,~~~\forall z\in S^{h}.
\end{eqnarray*}

\indent From Lemma 3.3 in \cite{yang5}, for any $(\psi,\varphi) \in
W^{m+1,p}(\Omega)\times W^{m-1,2}(\Omega)$, the following estimates
hold:
\begin{eqnarray}\label{s2.44r}
&&\||(\psi,\varphi)-I_{h}(\psi,\varphi)\||_{h}\leq C h^{m-1
+(\frac{1}{2}-\frac{1}{p})d}(\|\psi\|_{m+1,p,\Omega}+\|\varphi\|_{m-1,\Omega}),\\\label{s2.45r}
&&\|(\psi,\varphi)-I_{h}(\psi,\varphi)\|_{\mathbf{H}^{1}}\leq C
h^{m
+(\frac{1}{2}-\frac{1}{p})d}(\|\psi\|_{m+1,p,\Omega}+\|\varphi\|_{m-1,\Omega}).
\end{eqnarray}
From a Poincar$\acute{e}$-Friedrichs inequality \cite{brenner4} we get
\begin{eqnarray}\label{s2.46r}
&& \|(v,z)\|_{\mathbf{H}^{1}}= \|v\|_{1}+\|z\|_{-1}\lesssim \|v\|_{h}+\|z\|_{0}\nonumber\\
&&~~~\lesssim\|(v,z)\|_{h},~~~\forall
(v,z)\in \mathbf{H}_{h}.
\end{eqnarray}
Let $(v,z)=T_{h}(f,g)$ in (\ref{s2.33}), and we get
\begin{eqnarray}\label{s2.47r}
\|T_{h}(f,g)\|_{h}\leq C \|(f,g)\|_{\mathbf{H}^{1}},~~~\forall
(f,g)\in \mathbf{H}^{1}.
\end{eqnarray}

\indent{\bf Lemma 2.1.}~~ Let $(\psi,\varphi)$ and
$(\psi^{*},\varphi^{*})$ be the solution of (\ref{s2.9}) and
(\ref{s2.14}), respectively, and let $(\psi_{h},\varphi_{h})$ and
$(\psi_{h}^{*},\varphi_{h}^{*})$ be the $C^{0}$IPG approximation
solution of (\ref{s2.9}) and (\ref{s2.14}), respectively. Assume
that $(\psi,\varphi), (\psi^{*},\varphi^{*})\in
W^{m+1,p}(\Omega)\times H^{m-1}(\Omega)$ $(p\in \mathbb{S})$, then
\begin{eqnarray}\label{s2.48r}
&&\||(\psi,\varphi)-(\psi_{h},\varphi_{h})\||_{h}\lesssim
h^{m-1+(\frac{1}{2}-\frac{1}{p})d}(\|\psi\|_{m+1,p}+\|\varphi\|_{m-1}),\\\label{s2.49r}
&&\||(\psi^{*},\varphi^{*})-(\psi_{h}^{*},\varphi_{h}^{*})\||_{h}
\lesssim
h^{m-1+(\frac{1}{2}-\frac{1}{p})d}(\|\psi^{*}\|_{m+1,p}+\|\varphi^{*}\|_{m-1}),
\end{eqnarray}
furthermore assume $R(\Omega)$ holds, then
\begin{eqnarray}\label{s2.50r}
&&\|(\psi,\varphi)-(\psi_{h},\varphi_{h})\|_{\mathbf{H}^{1}}\lesssim
h^{m+(1-\frac{1}{p}-\frac{1}{p_{0}})d}(\|\psi\|_{m+1,p}+\|\varphi\|_{m-1}),
\\\label{s2.51r}
&&\|(\psi^{*},\varphi^{*})-(\psi_{h}^{*},\varphi_{h}^{*})\|_{\mathbf{H}^{1}}
\lesssim h^{m+(1-\frac{1}{p}-\frac{1}{p_{0}})d}
(\|\psi^{*}\|_{m+1,p}+\|\varphi^{*}\|_{m-1}).
\end{eqnarray}
\indent{\bf Proof.}~~ From (\ref{s2.31}), (\ref{s2.43r}),
(\ref{s2.29}) and (\ref{s2.44r}), we deduce
\begin{eqnarray*}
&&\|I_{h}(\psi,\varphi)-(\psi_{h},\varphi_{h})\|_{h}^{2} \lesssim
A_{h}(I_{h}(\psi,\varphi)-(\psi_{h},\varphi_{h}),
I_{h}(\psi,\varphi)-(\psi_{h},\varphi_{h}))\nonumber\\
&&~~~=A_{h}(I_{h}(\psi,\varphi)-(\psi,\varphi),
I_{h}(\psi,\varphi)-(\psi_{h},\varphi_{h}))\nonumber\\
&&~~~\lesssim  h^{m-1+d(\frac{1}{2}-\frac{1}{p})}
(\|\psi\|_{m+1,p}+\|\varphi\|_{m-1})
\|I_{h}(\psi,\varphi)-(\psi_{h},\varphi_{h})\|_{h},
\end{eqnarray*}
thus we get
\begin{eqnarray*}
&&\||(\psi,\varphi)-(\psi_{h},\varphi_{h})\||_{h}\leq
\||(\psi,\varphi)-I_{h}(\psi,\varphi)\||_{h}+\||I_{h}(\psi,\varphi)-(\psi_{h},\varphi_{h})\||_{h}\nonumber\\
&&~~~\lesssim
h^{m-1+d(\frac{1}{2}-\frac{1}{p})}(\|\psi\|_{m+1,p}+\|\varphi\|_{m-1}),
\end{eqnarray*}
which is the desired result (\ref{s2.48r}). By the same argument we
can
prove (\ref{s2.49r}).\\
Denote $e=\psi-\psi_{h}$. From (\ref{s2.39r}), (\ref{s2.42r}), (\ref{s2.43r}) with $z=0$,
(\ref{s2.29}), (\ref{s2.48r}), (\ref{s2.44r}) and (\ref{s2.40r}), we
deduce
\begin{eqnarray*}
&&|B'( e,  e)|= |a_{h}(Ke, e)|
= |a_{h}(e, Ke-I_{h}^{1}K e)|\nonumber\\
&&~~~\lesssim\||e\||_{h}\|| Ke-I_{h}^{1}K e\||_{h}\nonumber\\
&&~~~\lesssim
h^{m-1+(\frac{1}{2}-\frac{1}{p})d}\|\psi\|_{m+1,p} h^{1+(\frac{1}{2}-\frac{1}{p_{0}})d} \|K
e\|_{3,p_{0}}\nonumber\\
&&~~~\lesssim
h^{m+(1-\frac{1}{p}-\frac{1}{p_{0}})d}\|\psi\|_{m+1,p}\|e\|_{1},
\end{eqnarray*}
i.e.,
\begin{eqnarray}\label{s2.52r}
\|e\|_{1}\lesssim
h^{m+(1-\frac{1}{p}-\frac{1}{p_{0}})d}\|\psi\|_{m+1,p}.
\end{eqnarray}
From (\ref{s2.9}) and (\ref{s2.32}) we have
$\varphi=f \in H_{0}^{1}(\Omega)$ and
$$(\varphi-\varphi_{h}, z)_{0}=0, ~~~\forall z\in S^{h}.$$
So
\begin{eqnarray}\label{s2.53r}
\|\varphi-\varphi_{h}\|_{-1}\lesssim
h^{m}\|\varphi\|_{m-1}.
\end{eqnarray}
From (\ref{s2.52r}) and (\ref{s2.53r}) we get the desired result (\ref{s2.50r}). By the same argument
we can prove (\ref{s2.51r}).
The proof is completed.~~~$\Box$\\

Based on Lemma 2.1, using argument as  Theorem 3.3 and  Theorem 3,4 in \cite{yang5} we can prove the following a priori error estimates for the eigenvalue problem.\\

\indent{\bf Theorem 2.1.}~ Assume that $R(\Omega)$ holds and $n\in W^{1,\infty}(\Omega)\cap H^{2}(\Omega)$, then
\begin{eqnarray}
\label{s2.38}
&&|(\frac{1}{q}\sum\limits_{j=i}^{i+q-1}\lambda_{j,h}^{-1})^{-1}-\lambda|\lesssim
\|(T-T_{h})|_{ran(E)}\|_{\mathbf{H}^{1}};
\end{eqnarray}
assume $ran(E)\subset W^{m+1,p}(\Omega)\times H^{m-1}(\Omega)$
$(p\in \mathbb{S})$, then
\begin{eqnarray}\label{s2.39}
\|(T-T_{h})|_{ran(E)}\|_{\mathbf{H}^{1}}\lesssim
h^{m+(1-\frac{1}{p}-\frac{1}{p_{0}})d};
\end{eqnarray}
furthermore, assume that $(u_{h},\omega_{h})$ is an eigenfunction
corresponding to $\lambda_{h}$ and $\|(u_{h},\omega_{h})\|_h=1$,
then there exists eigenfunction $(u,\omega)$ corresponding to
$\lambda$ such that
\begin{eqnarray}\label{s2.40}
&&\|(u_{h},\omega_{h})-(u,\omega)\|_{\mathbf{H}^{1}} \lesssim
h^{\frac{m}{\alpha}+(1-\frac{1}{p}-\frac{1}{p_{0}})\frac{d}{\alpha}},\\\label{s2.41}
&& \||(u_{h},\omega_{h})-(u,\omega)\||_{h}\lesssim
h^{\frac{m-1}{\alpha}+(\frac{1}{2}-\frac{1}{p})\frac{d}{\alpha}}.
\end{eqnarray}

\indent  In addition, when the eigenvalue $\lambda$ is non-defective, for
$(u^{*},\omega^{*})\in ran(E^{*})$ with
$\|(u^{*},\omega^{*})\|_h=1$, there exists
$(u_{h}^{*},\omega_{h}^{*})\in ran(E_{h}^{*})$ such that
\begin{eqnarray}\label{s2.42}
&&\|(u^{*},\omega^{*})-(u_{h}^{*},\omega_{h}^{*})\|_{\mathbf{H}^{1}}\lesssim
h^{m+(1-\frac{1}{p}-\frac{1}{p_{0}})d},\\\label{s2.43}
 &&\||(u^{*},\omega^{*})-(u_{h}^{*},\omega_{h}^{*})\||_{h}\lesssim
h^{m-1+(\frac{1}{2}-\frac{1}{p})d};
\end{eqnarray}
for $(u_{h}^{*},\omega_{h}^{*})\in ran(E_{h}^{*})$ with
$\|(u_{h}^{*},\omega_{h}^{*})\|_h=1$, there exists
$(u^{*},\omega^{*})\in ran(E^{*})$ such that
\begin{eqnarray}\label{s2.44}
&&\|(u_{h}^{*},\omega_{h}^{*})-(u^{*},\omega^{*})\|_{\mathbf{H}^{1}}\lesssim
h^{m+(1-\frac{1}{p}-\frac{1}{p_{0}})d},\\\label{s2.45}
 &&\||(u_{h}^{*},\omega_{h}^{*})-u^{*},\omega^{*})\||_{h}\lesssim
h^{m-1+(\frac{1}{2}-\frac{1}{p})d},\\\label{s2.46}
&&|\lambda_{h}-\lambda|\lesssim
h^{2m-2+2(\frac{1}{2}-\frac{1}{p})d}.
\end{eqnarray}

\section{The a posteriori error analysis of $C^{0}$IPG discrete scheme for  the source problem (\ref{s2.9}) }

\indent In 2012, Brenner \cite{brenner1} proposed and analyzed the a
posteriori error estimates of $C^0$IPG methods for biharmonic
equation. Based on \cite{brenner1}, in this section we discuss a posteriori error estimates
of $C^0$IPG discrete scheme (\ref{s2.32}) for the source problem (\ref{s2.9}).\\
\indent Denote $$F=F(f,g)=-\Delta(\frac{1}{n-1}f)
-\frac{n}{n-1}\Delta f
-\frac{n}{n-1}g,~~~in~\kappa,$$
where $f,g\in W^{3,p}(\Omega, \pi_{h})$,
and denote
\begin{eqnarray}\label{s3.1}
&&\eta_{\kappa}(F,\psi_{h})=h^2_{\kappa}\| F-\Delta(\frac{1}{n-1}\Delta\psi_{h})\|_{0,\kappa},
~\forall \kappa\in\mathcal{\pi}_h,\\\label{s3.2}
&&\eta_{\ell,1}(\psi_{h})=\frac{\sigma}{\hat{\ell}^{\frac{1}{2}}}\left\|\left[\left[\frac{\partial
\psi_h}{\partial\gamma_{\ell}}\right]\right]\right\|_{0,\ell},
~\forall \ell\in\mathcal{E},\\\label{s3.3}
&&\eta_{\ell,2}(\psi_{h})=\mu\hat{\ell}^{\frac{1}{2}}\left\|\left[\left[\frac{\partial^2\psi_h}{\partial\gamma^2_{\ell}}
\right]\right]\right\|_{0,\ell},~~
~\forall\ell\in\mathcal{E}^{i},\\
\label{s3.4}
&&
\eta_{\ell,3}(\psi_{h})=\hat{\ell}^{\frac{3}{2}}\left\|\left[\left[\frac{\partial(\frac{1}{n-1}\Delta
\psi_h)}{\partial\gamma_{\ell}}\right]\right]\right\|_{0,\ell},
~~~\forall\ell\in\mathcal{E}^{i},\\\label{s3.5}
&&\eta_{\ell,4}(\psi_{h})=\hat{\ell}^{\frac{1}{2}}\left\|\left[\left[(\frac{1}{n-1}-\mu)\Delta \psi_h\right]\right]\right\|_{0,\ell},~~
~\forall\ell\in\mathcal{E}^{i}.
\end{eqnarray}
Then the residual-based error indicator $\eta_h$ is defined by
\begin{eqnarray}\label{s3.6}
&&\eta^2_h(F,\psi_{h},\varphi_{h},\kappa)=\eta^2_{\kappa}(F,\psi_{h})
+\sum_{\ell\in\mathcal{E}^{b}\cap\partial\kappa}\eta^2_{\ell,1}(\psi_{h})\nonumber\\
&&~~~
+\frac{1}{2}\sum_{\ell\in\mathcal{E}^{i}\cap\kappa}
\{\eta^2_{\ell,1}(\psi_{h})+\eta^2_{\ell,2}(\psi_{h})+\eta^2_{\ell,3}(\psi_{h})+\eta^2_{\ell,4}(\psi_{h})\}\nonumber\\
&&~~~+\sum\limits_{\ell\in
\mathcal{E}\cap\partial\kappa}\|\frac{n+1}{n-1}\|_{0,\ell}h_{\ell}^{4}\eta_{\ell,1}^2(f)+\|f-\varphi_{h}\|_{0,\kappa}^{2}+h^4\sum\limits_{\ell\in
\mathcal{E}\cap\partial\kappa}\eta^2_{\ell,1}(\varphi_{h}),\\\label{s3.7}
&&\eta^2_h(F,\psi_{h},\varphi_{h},\Omega)=\sum_{\kappa\in\pi_{h}}\eta^2_h(F,\psi_{h},\varphi_{h},\kappa).
\end{eqnarray}

\indent Let $P_{j}(\Omega,\pi_{h})$ be the space of piecewise polynomial functions of degree $\leq j$ and $\widetilde{g}\in P_{j}(\Omega,\pi_{h})$
denote the $L^{2}$ orthogonal projection of $g$. And denote
\begin{eqnarray}\label{s3.8}
&&\widehat{F}=-\Delta(\widetilde{\frac{1}{n-1}f})
-\widetilde{\frac{n}{n-1}\Delta f}
-\widetilde{\frac{n}{n-1}g},\\\label{s3.9}
&&\widehat{\eta}_{\kappa}(F,\psi_{h})=h^2_{\kappa}\|\widehat{F}-\Delta(\widetilde{\frac{1}{n-1}}\Delta\psi_{h})\|_{0,\kappa},
~\forall \kappa\in\mathcal{\pi}_h,\\\label{s3.10}
&&
\widehat{\eta}_{\ell,3}(\psi_{h})=\hat{\ell}^{\frac{3}{2}}\left\|\left[\left[\frac{\partial(\widetilde{\frac{1}{n-1}}\Delta
\psi_h)}{\partial\gamma_{\ell}}\right]\right]\right\|_{0,\ell},
~~~\forall\ell\in\mathcal{E}^{i},\\\label{s3.11}
&&\widehat{\eta}_{\ell,4}(\psi_{h})=\hat{\ell}^{\frac{1}{2}}\left\|\left[\left[(\widetilde{\frac{1}{n-1}-\mu})\Delta \psi_h\right]\right]\right\|_{0,\ell},~~
~\forall\ell\in\mathcal{E}^{i}.
\end{eqnarray}

The data oscillations are defined by
\begin{eqnarray}\label{s3.12}
&&Osc_{j}(F)=(\sum\limits_{\kappa\in \pi_{h}}h_{\kappa}^4\|F-\widehat{F}\|^2_{0,\kappa})^{\frac{1}{2}},\\\label{s3.13}
&&Osc_{j}(\eta_{\ell,3})=(
\sum\limits_{\kappa\in \pi_{h}}\sum\limits_{\ell\in
\mathcal{E}^i\cap\partial\kappa} (\eta_{\ell,3}(\psi_{h})-\widehat{\eta}_{\ell,3}(\psi_{h}))^2)^{\frac{1}{2}},\\\label{s3.14}
&&Osc_{j}(\eta_{\ell,4})=(
\sum\limits_{\kappa\in \pi_{h}}\sum\limits_{\ell\in
\mathcal{E}^i\cap\partial\kappa} (\eta_{\ell,4}(\psi_{h})-\widehat{\eta}_{\ell,4}(\psi_{h}))^2)^{\frac{1}{2}}.
\end{eqnarray}

\indent{\bf Theorem 3.1.}~~ Let $(\psi,\varphi)$ and
$(\psi_{h},\varphi_{h})$ be the solution of (\ref{s2.9}) and
(\ref{s2.32}), respectively.  Assume
that $R(\Omega)$ holds and $n\in W^{1,\infty}(\Omega)\cap H^{2}(\Omega)$, then
\begin{eqnarray}\label{s3.15}
\|(\psi,\varphi)-(\psi_{h},\varphi_{h})\|_{h}\lesssim
\eta_h(F,\psi_{h},\varphi_{h},\Omega).
\end{eqnarray}
\indent{\bf Proof.}~~
Brenner introduced the enriching operator $E_{h}:
S^{h}\to H^{2}(\Omega)$ and proved (see (4.4) in \cite{brenner1})
\begin{eqnarray}\label{s3.16}
&&~~~\sum\limits_{\kappa\in\pi_{h}}(h_{\kappa}^{-4}\|v-E_{h}v\|_{0,\kappa}^{2}+h_{\kappa}^{-2}|v-E_{h}v|_{1,\kappa}^{2}+|v-E_{h}v|_{2,\kappa}^{2})\nonumber\\
&&\lesssim
\sum\limits_{\kappa\in\pi_{h}}\frac{1}{\hat{\ell}}\|[[\partial
v/\partial\gamma_{\ell}]]\|_{0,\ell}^{2},~~~\forall v\in S^{h}.
\end{eqnarray}

\indent Denote $\mathbf{E}_{h}(u_{h},\omega_{h})=(E_{h}u_{h}, E_{h}\omega_{h})$.\\
\indent Due to (\ref{s2.26}) we need to bound $\sigma\sum\limits_{\ell\in
\mathcal {E}}\frac{1}{\hat{\ell}} \|[[\frac{\partial
(\psi-\psi_{h})}{\partial \gamma_{\ell}}]]\|_{0,\ell}^{2}$ and
$\sum\limits_{\kappa\in\pi_{h}} \|
\psi-\psi_{h}\|_{2,\kappa}^{2}+
\sum\limits_{\kappa\in\pi_{h}}\|\varphi-\varphi_{h}\|_{0,\kappa}^{2}$.\\
\indent Since $\sigma>1$, from (\ref{s3.2}) we get
\begin{eqnarray}\label{s3.17}
\sigma\sum\limits_{\ell\in \mathcal {E}}\frac{1}{\hat{\ell}}
\|[[\frac{\partial (\psi-\psi_{h})}{\partial
\gamma_{\ell}}]]\|_{0,\ell}^{2}=\sigma\sum\limits_{\ell\in \mathcal {E}}\frac{1}{\hat{\ell}}
\|[[\frac{\partial \psi_{h}}{\partial
\gamma_{\ell}}]]\|_{0,\ell}^{2}\leq \sum\limits_{\ell\in \mathcal
{E}}\eta_{\ell,1}^{2}.
\end{eqnarray}
From (\ref{s3.16}) and (\ref{s3.2}) we have
\begin{eqnarray}\label{s3.18}
&&~~~\sum\limits_{\kappa\in\pi_{h}} \| \psi-\psi_{h}\|_{2,\kappa}^{2} +
\sum\limits_{\kappa\in\pi_{h}}\|\varphi-\varphi_{h}\|_{0,\kappa}^{2}\nonumber\\
&&\leq 2\sum\limits_{\kappa\in\pi_{h}} (\|
\psi-E_{h}\psi_{h}\|_{2,\kappa}^{2} +|\psi_{h}-E_{h}\psi_{h}\|_{2,\kappa}^{2})\nonumber\\
&&~~~+2\sum\limits_{\kappa\in\pi_{h}}(\|\varphi-E_{h}\varphi_{h}\|_{0,\kappa}^{2}+\|\varphi_{h}-E_{h}\varphi_{h}\|_{0,\kappa}^{2})\nonumber\\
&&\lesssim\|(\psi,\varphi)-\mathbf{E}_{h}(\psi_{h},\varphi_{h})\|_{\mathbf{H}}^{2}+\sum\limits_{\ell\in
\mathcal
{E}}\frac{1}{\hat{\ell}}\eta_{\ell,1}^{2}(\psi_{h})+h^{4}\sum\limits_{\ell\in
\mathcal {E}}\eta_{\ell,1}^{2}(\varphi_{h}).
\end{eqnarray}
By duality we have
\begin{eqnarray}\label{s3.19}
\|(\psi,\varphi)-\mathbf{E}_{h}(\psi_{h},\varphi_{h})\|_{\mathbf{H}}\thickapprox\sup\limits_{(v,z)\in
\mathbf{H}\setminus
\{0\}}\frac{A((\psi,\varphi)-\mathbf{E}_{h}(\psi_{h},\varphi_{h}),
(v,z))}{\|(v,z)\|_{\mathbf{H}}}.
\end{eqnarray}
Denote
\begin{eqnarray}\label{s3.20}
&&A_{\kappa}((u,\omega),(v,z))=\int\limits_{\kappa}(\frac{1}{n-1}-\mu)\Delta u, \Delta
\overline{v}dx\nonumber\\
&&~~~+\mu\int\limits_{\kappa}D^{2}u:D^{2}\bar{v}dx+\int\limits_{\kappa}\omega \overline{z}dx.
\end{eqnarray}
From (\ref{s2.6}), (\ref{s3.20}), (\ref{s2.9}) and  (\ref{s2.32}) we get
\begin{eqnarray*}
&&~~~A((\psi,\varphi)-\mathbf{E}_{h}(\psi_{h},\varphi_{h}), (v,z))\nonumber\\
&&=
\sum\limits_{\kappa\in\pi_{h}}A_{\kappa}((\psi_{h},\varphi_{h})-\mathbf{E}_{h}(\psi_{h},\varphi_{h}),
(v,z))-\sum\limits_{\kappa\in\pi_{h}}A_{\kappa}((\psi_{h},\varphi_{h}), (v,z)-I_{h}(v,z))\nonumber
\end{eqnarray*}

\begin{eqnarray}\label{s3.21}
&&~~~+A((\psi,\varphi),
(v,z))-\sum\limits_{\kappa\in\pi_{h}}A_{\kappa}((\psi_{h},\varphi_{h}),
I_{h}(v,z))\nonumber\\
&&=
\sum\limits_{\kappa\in\pi_{h}}A_{\kappa}((\psi_{h},\varphi_{h})-\mathbf{E}_{h}(\psi_{h},\varphi_{h}),
(v,z))\nonumber\\
&&~~~-\sum\limits_{\kappa\in\pi_{h}}A_{\kappa}((\psi_{h},\varphi_{h}), (v,z)-I_{h}(v,z))\nonumber\\
&&~~~+A_{h}((\psi_{h},\varphi_{h}),
I_{h}(v,z))-\sum\limits_{\kappa\in\pi_{h}}A_{\kappa}((\psi_{h},\varphi_{h}),
I_{h}(v,z))\nonumber\\
&&~~~+B((f,g),(v,z)-I_{h}(v,z))\equiv
I_{1}-I_{2}+I_{3}-I_{4}+I_{5}.
\end{eqnarray}
We have
\begin{eqnarray}\label{s3.22}
&&I_{2}=\sum\limits_{\kappa\in\pi_{h}}A_{\kappa}((\psi_{h},\varphi_{h}),
(v,z)-I_{h}(v,z))\nonumber\\
&&~~~=\sum\limits_{\kappa\in\pi_{h}} \int\limits_{\kappa}
(\frac{1}{n-1}-\mu)\Delta \psi_{h} \Delta \overline{(v-I_{h}^{1}v)}dx+
\mu\int\limits_{\kappa} D^{2}\psi_{h}: D^{2}
\overline{(v-I_{h}^{1}v)}dx
\nonumber\\
&&~~~~~~+\sum\limits_{\kappa\in\pi_{h}}\int\limits_{\kappa}\varphi_{h}
\overline{(z-I_{h}^{2}z)}dx\equiv J_{1}+J_{2}+J_{3},
\end{eqnarray}
and by the Green's formula we have
\begin{eqnarray}\label{s3.23}
&&J_{1}= \sum\limits_{\kappa\in\pi_{h}}-\int\limits_{\kappa}
\nabla[(\frac{1}{n-1}-\mu)\Delta\psi_{h}]
\overline{\nabla(v-I_{h}^{1}v)}dx\nonumber\\
&&~~~~~~+\sum\limits_{\kappa\in\pi_{h}}\int_{\partial\kappa}(\frac{1}{n-1}-\mu)\Delta\psi_{h}\overline{\frac{\partial(v-I_{h}^{1}v)}{\partial\gamma}}
ds\nonumber\\
&&~~~=\sum\limits_{\kappa\in\pi_{h}}\int\limits_{\kappa}
\Delta[(\frac{1}{n-1}-\mu)\Delta\psi_{h}]
\overline{(v-I_{h}^{1}v)}dx\nonumber\\
&&~~~~~~-\sum\limits_{\kappa\in\pi_{h}}\int_{\partial\kappa}\nabla[(\frac{1}{n-1}-\mu)\Delta\psi_{h}]\overline{(v-I_{h}^{1}v)\cdot\gamma}ds\nonumber\\
&&~~~~~~+\sum\limits_{\kappa\in\pi_{h}}\int_{\partial\kappa}(\frac{1}{n-1}-\mu)
\Delta\psi_{h}\overline{\frac{\partial(v-I_{h}^{1}v)}{\partial\gamma}}ds\nonumber\\
&&~~~= \sum\limits_{\kappa\in\pi_{h}}\int\limits_{\kappa}
\Delta[(\frac{1}{n-1}-\mu)\Delta\psi_{h}]
\overline{(v-I_{h}^{1}v)}dx\nonumber\\
&&~~~~~~+\sum\limits_{\ell\in\mathcal{E}}\int_{\ell}[[\nabla[(\frac{1}{n-1}-\mu)\Delta\psi_{h}]\cdot\gamma]]\overline{(v-I_{h}^{1}v)}ds\nonumber\\
&&~~~~~~-\sum\limits_{\ell\in\mathcal{E}}\int_{\ell}\{\{(\frac{1}{n-1}-\mu)\Delta\psi_{h}\}\}\overline{[[\frac{\partial(v-I_{h}^{1}v)}
{\partial\gamma_{\ell}}]]}ds\nonumber\\
&&~~~~~~-\sum\limits_{\ell\in\mathcal
{E}^i}\int_{\ell}[[(\frac{1}{n-1}-\mu)\Delta\psi_{h}]]\overline{\{\{\frac{\partial(v-I_{h}^{1}v)}
{\partial\gamma_{\ell}}\}\}}ds,
\end{eqnarray}
and by the Green's formula (see also (7.10) in \cite{brenner1}) we have
\begin{eqnarray}\label{s3.24}
&&J_{2}= \mu\{\sum\limits_{\kappa\in\pi_{h}}\int\limits_{\kappa}
(\Delta^{2}\psi_{h}) \overline{(v-I_{h}^{1}v)}dx +\sum\limits_{\ell\in
\mathcal {E}^i}\int_{\ell} [[\frac{\partial
(\Delta\psi_{h})}{\partial\gamma_{\ell}}]]
\overline{(v-I_{h}^{1}v)}ds\nonumber\\
&&~~~+\sum\limits_{\ell\in \mathcal {E}}\int_{\ell}
\{\{\frac{\partial^{2}\psi_{h}}{\partial\gamma_{\ell}^{2}}\}\}
\overline{[[\frac{\partial I_{h}^{1}v}{\partial\gamma_{\ell}}]]}ds
 -\sum\limits_{\ell\in \mathcal {E}^i}\int_{\ell}
[[\frac{\partial^{2}\psi_{h}}{\partial\gamma_{\ell}^{2}}]]
\overline{\{\{\frac{\partial(v-I_{h}^{1}v)}{\partial\gamma_{\ell}}\}\}}ds
\nonumber\\
&&~~~ - \sum\limits_{\ell\in \mathcal {E}^i}\int_{\ell}
[[\frac{\partial^{2}\psi_{h}}{\partial\gamma_{\ell}\partial
t_{\ell}}]] \frac{\partial(\overline{v-I_{h}^{1}v})}{\partial t_{\ell}}ds\},
\end{eqnarray}
By (\ref{s2.24}) we get
\begin{eqnarray}\label{s3.25}
&&I_{3}-I_{4}=
 \sum\limits_{\ell\in \mathcal {E}}\int\limits_{\ell}
\{\{(\frac{1}{n-1}-\mu)\Delta \psi_{h}\}\}[[\frac{\partial
\overline{I_{h}^{1}v}}{\partial
\gamma_{\ell}}]]\nonumber\\
&&~~~~~~+\{\{(\frac{1}{n-1}-\mu)\Delta
I_{h}^{1}v\}\}[[\frac{\partial
\overline{\psi_{h}}}{\partial \gamma_{\ell}}]]ds\nonumber\\
&&~~~~~~+ \mu\sum\limits_{\ell\in \mathcal {E}}\int\limits_{\ell}
\{\{\frac{\partial^{2}\psi_{h}}{\partial\gamma_{\ell}^{2}}\}\}[[\frac{\partial
\overline{I_{h}^{1}v}}{\partial
\gamma_{\ell}}]]+\{\{\frac{\partial^{2}I_{h}^{1}v}{\partial\gamma_{\ell}^{2}}\}\}[[\frac{\partial
\overline{\psi_{h}}}{\partial \gamma_{\ell}}]]ds\nonumber\\
&&~~~~~~ + \sigma\sum\limits_{\ell\in \mathcal
{E}}\frac{1}{\hat{\ell}}\int\limits_{\ell}[[\frac{\partial
\psi_{h}}{\partial \gamma_{\ell}}]] [[\frac{\partial
\overline{I_{h}^{1}v}}{\partial \gamma_{\ell}}]]ds,
\end{eqnarray}
from  the Green's formula we get
\begin{eqnarray}\label{s3.26}
&&~~~B((f,g),(v,z))\nonumber\\
&&=(\nabla(\frac{1}{n-1}f), \nabla v)_{0}+(\nabla f,
\nabla(\frac{n}{n-1}v))_{0}-(g, \frac{n}{n-1}v
)_{0}+(f,z)_{0}\nonumber\\
&&=-\sum\limits_{\kappa}\int\limits_{\kappa}\Delta(\frac{1}{n-1}f)
\overline{v}dx-\sum\limits_{\kappa}\int\limits_{\kappa}\frac{n}{n-1}\Delta f
\overline{v}dx-(\frac{n}{n-1}g, v
)_{0}+(f,z)_{0}\nonumber\\
&&~~~+
\sum\limits_{\kappa}\int\limits_{\partial\kappa}\frac{\partial(\frac{1}{n-1}f)}{\partial\gamma}
\overline{v}ds+\sum\limits_{\kappa}\int\limits_{\partial\kappa}\frac{n}{n-1}\frac{\partial
f}{\partial\gamma}\overline{ v}ds\nonumber\\
&&\equiv \sum\limits_{\kappa}\int\limits_{\kappa}F\overline{v}dx
+(f,z)_{0}+\sum\limits_{\ell\in
\mathcal{E}}\int\limits_{\ell}[[\frac{\partial(\frac{1}{n-1}f)}{\partial\gamma}]]
\overline{v}ds\nonumber\\
&&~~~+\sum\limits_{\ell\in
\mathcal{E}}\int\limits_{\ell}[[\frac{n}{n-1}\frac{\partial
f}{\partial\gamma}]] \overline{v}ds,
\end{eqnarray}
thus
\begin{eqnarray}\label{s3.27}
&&I_{5}=B((f,g),(v,z)-I_{h}(v,z))
= \sum\limits_{\kappa}\int\limits_{\kappa}F(\overline{v-I_{h}^{1}v})dx\nonumber\\
&&~~~~~~+(f,z-I_{h}^{2}z)_{0}
+\sum\limits_{\ell\in
\mathcal{E}}\int\limits_{\ell}[[\frac{n+1}{n-1}\frac{\partial
f}{\partial\gamma}]] (\overline{v-I_{h}^{1}v})ds.
\end{eqnarray}
Substituting (\ref{s3.22}),(\ref{s3.25}) and (\ref{s3.27}) into
(\ref{s3.21}), we obtain
\begin{eqnarray}\label{s3.28}
&&~~~A((\psi,\varphi)-\mathbf{E}_{h}(\psi_{h},\varphi_{h}), (v,z))\nonumber\\
&&= I_{1}+\sum\limits_{\kappa\in\pi_{h}}\int\limits_{\kappa}
(F-\Delta(\frac{1}{n-1}\Delta\psi_{h}))\overline{(v-I_{h}^{1}v)}+(f-\varphi_{h})\overline{(z-I_{h}^{2}z)}dx\nonumber\\
&&~~~-\sum\limits_{\ell\in\mathcal{E}}\int_{\ell}[[\nabla[(\frac{1}{n-1}-\mu)\Delta\psi_{h}]\cdot\gamma]]\overline{(v-I_{h}^{1}v)}ds\nonumber\\
&&~~~+\sum\limits_{\ell\in\mathcal{E}}\int_{\ell}\{\{(\frac{1}{n-1}-\mu)\Delta\psi_{h}\}\}\overline{[[\frac{\partial(v-I_{h}^{1}v)}
{\partial\gamma_{\ell}}]]}ds\nonumber\\
&&~~~+\sum\limits_{\ell\in\mathcal
{E}^i}\int_{\ell}[[(\frac{1}{n-1}-\mu)\Delta\psi_{h}]]\overline{\{\{\frac{\partial(v-I_{h}^{1}v)}
{\partial\gamma_{\ell}}\}\}}ds\nonumber\\
&&~~~-\mu\sum\limits_{\ell\in \mathcal {E}^i}\int_{\ell}
[[\frac{\partial (\Delta\psi_{h})}{\partial\gamma_{\ell}}]]
\overline{(v-I_{h}^{1}v)}ds-\mu\sum\limits_{\ell\in \mathcal {E}}\int_{\ell}
\{\{\frac{\partial^{2}\psi_{h}}{\partial\gamma_{\ell}^{2}}\}\}
\overline{[[\frac{\partial I_{h}^{1}v}{\partial\gamma_{\ell}}]]}ds\nonumber\\
&&~~~ +\mu
\sum\limits_{\ell\in \mathcal {E}^i}\int_{\ell}
[[\frac{\partial^{2}\psi_{h}}{\partial\gamma_{\ell}^{2}}]]
\overline{\{\{\frac{\partial(v-I_{h}^{1}v)}{\partial\gamma_{\ell}}\}\}}ds
+ \mu\sum\limits_{\ell\in \mathcal {E}^i}\int_{\ell}
[[\frac{\partial^{2}\psi_{h}}{\partial\gamma_{\ell}\partial
t_{\ell}}]] \frac{\partial(\overline{v-I_{h}^{1}v})}{\partial
t_{\ell}}ds\nonumber\\
&&~~~+\sum\limits_{\ell\in \mathcal {E}}\int\limits_{\ell}
\{\{(\frac{1}{n-1}-\mu)\Delta \psi_{h}\}\}[[\frac{\partial
\overline{I_{h}^{1}v}}{\partial
\gamma_{\ell}}]]ds\nonumber\\
&&~~~+\sum\limits_{\ell\in \mathcal {E}}\int\limits_{\ell}\{\{(\frac{1}{n-1}-\mu)\Delta
I_{h}^{1}v\}\}[[\frac{\partial
\overline{\psi_{h}}}{\partial \gamma_{\ell}}]]ds+ \mu\sum\limits_{\ell\in \mathcal {E}}\int\limits_{\ell}
\{\{\frac{\partial^{2}\psi_{h}}{\partial\gamma_{\ell}^{2}}\}\}[[\frac{\partial
\overline{I_{h}^{1}v}}{\partial
\gamma_{\ell}}]]ds\nonumber\\
&&~~~+\mu\sum\limits_{\ell\in \mathcal {E}}\int\limits_{\ell}\{\{\frac{\partial^{2}I_{h}^{1}v}{\partial\gamma_{\ell}^{2}}\}\}[[\frac{\partial
\overline{\psi_{h}}}{\partial \gamma_{\ell}}]]ds + \sigma\sum\limits_{\ell\in \mathcal
{E}}\frac{1}{\hat{\ell}}\int\limits_{\ell}[[\frac{\partial
\psi_{h}}{\partial \gamma_{\ell}}]] [[\frac{\partial
\overline{I_{h}^{1}v}}{\partial \gamma_{\ell}}]]ds\nonumber\\
&&~~~+
\sum\limits_{\ell\in
\mathcal{E}}\int\limits_{\ell}[[\frac{n+1}{n-1}\frac{\partial
f}{\partial\gamma}]] (\overline{v-I_{h}^{1}v})ds\nonumber\\
&&\equiv I_{1}+G_{2}+G_{3}+\cdots+G_{15}.
\end{eqnarray}
By (\ref{s3.20}), the Schwarz inequality, (\ref{s3.16}) and (\ref{s3.2}) we get
\begin{eqnarray}\label{s3.29}
&&|I_{1}|=|\sum\limits_{\kappa\in\pi_{h}}A_{\kappa}((\psi_{h},\varphi_{h})-\mathbf{E}_{h}(\psi_{h},\varphi_{h}),
(v,z))|\nonumber\\
&&~~~\lesssim\sum\limits_{\kappa\in\pi_{h}}(\|\frac{1}{n-1}\|_{0,\infty,\kappa}|\psi_{h}-E_{h}\psi_{h}|_{2,\kappa}+\|\varphi_{h}-E_{h}\varphi_{h}\|_{0,\kappa})
\|(v,z)\|_{\mathbf{H}}\nonumber\\
&&~~~\lesssim \sum\limits_{\ell\in
\mathcal{E}}\int\limits_{\ell}(
\|\frac{1}{n-1}\|_{0,\infty,\kappa}^2\eta_{\ell,1}(\psi_{h})^{2}+\eta_{\ell,1}(\varphi_{h})^{2})^{\frac{1}{2}}\|(v,z)\|_{\mathbf{H}},
\end{eqnarray}
by (\ref{s3.1}) we get
\begin{eqnarray*}
&&|G_{2}|\lesssim (\sum\limits_{\kappa\in
\pi_{h}}h_{\kappa}^4\|  F-\Delta(\frac{1}{n-1}\Delta\psi_{h}) \|_{0,\kappa}^{2})^{\frac{1}{2}}
(\sum\limits_{\kappa\in
\pi_{h}}h_{\kappa}^{-4}\| v-I_{h}^{1}v  \|_{0,\kappa}^{2})^{\frac{1}{2}}\nonumber\\
&&~~~~~~+\|f-\varphi_{h}\|_{0}\|z-I_{h}^{2}z\|_{0}\lesssim (\sum\limits_{\kappa\in
\pi_{h}}\eta_{\kappa}^{2})^{\frac{1}{2}}|v|_{2}+\|f-\varphi_{h}\|_{0}\|z\|_{0},
\end{eqnarray*}
by (\ref{s3.4}) we get
\begin{eqnarray*}
&&|G_{3}+G_{6}|\lesssim (\sum\limits_{\ell\in
\mathcal{E}^i}\hat{\ell}^3\|[[\nabla(\frac{1}{n-1}\Delta\psi_{h})\cdot\gamma]]\|_{0,\ell}^{2})^{\frac{1}{2}}
(\sum\limits_{\ell\in
\mathcal{E}^i}\hat{\ell}^{-3}\|v-I_{h}^{1}v\|_{0,\ell}^{2})^{\frac{1}{2}}\nonumber\\
&&~~~\lesssim (\sum\limits_{\ell\in
\mathcal{E}^i}\eta_{\ell,3}^{2})^{\frac{1}{2}}|v|_{2},
\end{eqnarray*}
we see
\begin{eqnarray*}
G_{4}+G_{10}=0,~~~G_{7}+G_{12}=0,
\end{eqnarray*}
by (\ref{s3.5}) we get
\begin{eqnarray*}
|G_{5}|\lesssim (\sum\limits_{\ell\in
\mathcal{E}^i}h_{\ell}[[(\frac{1}{n-1}-\mu)\Delta\psi_{h}]]^{2})^{\frac{1}{2}}|v|_{2}
\lesssim (\sum\limits_{\ell\in
\mathcal{E}^i}\eta_{\ell£¬4}^{2})^{\frac{1}{2}}|v|_{2},
\end{eqnarray*}
by (\ref{s3.3}) we get
\begin{eqnarray*}
|G_{8}|\lesssim (\sum\limits_{\ell\in
\mathcal{E}^i}\eta_{\ell,2}(\psi_{h})^{2})^{\frac{1}{2}}|v|_{2},
\end{eqnarray*}
by (\ref{s3.2}), the trace theorem with scaling and a standard inverse estimate, we deduce
\begin{eqnarray*}
&&|G_{9}|\lesssim \mu(\sum\limits_{\ell\in
\mathcal{E}^i}\hat{\ell}\|[[\frac{\partial^2\psi_{h}}{\partial\gamma_{\ell}\partial t_{\ell}}]]\|_{0,\ell}^2)^{\frac{1}{2}}
(\sum\limits_{\ell\in
\mathcal{E}^i}\hat{\ell}^{-1}\|\frac{\partial(v-I_{h}^{1}v)}{\partial t_{\ell}}\|_{0,\ell}^{2})^{\frac{1}{2}}\\
&&~~~\lesssim \mu(\sum\limits_{\ell\in
\mathcal{E}^i}\eta_{\ell,1}(\psi_{h})^{2})^{\frac{1}{2}}|v|_{2},
\end{eqnarray*}
by (\ref{s3.2}) we get
\begin{eqnarray*}
&&|G_{11}|\lesssim \frac{1}{\sigma}(\sum\limits_{\ell\in
\mathcal{E}}\|\frac{1}{n-1}-\mu\|_{0,\infty,\ell}\eta_{\ell,1}(\psi_{h})^{2})^{\frac{1}{2}}|v|_{2},\nonumber\\
&&|G_{13}|\lesssim \frac{\mu}{\sigma}(\sum\limits_{\ell\in
\mathcal{E}}\eta_{\ell,1}(\psi_{h})^{2})^{\frac{1}{2}}|v|_{2},\nonumber\\
&&|G_{14}|\lesssim (\sum\limits_{\ell\in
\mathcal{E}}\eta_{\ell,1}(\psi_{h})^{2})^{\frac{1}{2}}|v|_{2},\nonumber\\
&&|G_{15}|\lesssim (\sum\limits_{\ell\in
\mathcal{E}}\|\frac{n+1}{n-1}\|_{0,\infty,\ell}h_{\ell}^{4}\eta_{\ell,1}(f)^{2})^{\frac{1}{2}}|v|_{2}.\nonumber
\end{eqnarray*}
Substituting these estimates into
(\ref{s3.28}), we obtain
\begin{eqnarray}\label{s3.30}
A((\psi,\varphi)-\mathbf{E}_{h}(\psi_{h},\varphi_{h}), (v,z))\lesssim \eta_h(F,\psi_{h},\varphi_{h},\Omega)\|(v,z)\|_{\mathbf{H}}.
\end{eqnarray}
\indent Combining (\ref{s3.17})-(\ref{s3.19}) and (\ref{s3.30}) we obtain (\ref{s3.15}).~~~$\Box$\\

Using the argument in Theorem 8 in \cite{brenner1}, we can prove the following theorem.\\
\indent{\bf Theorem 3.2.}~~ Under the condition of Theorem 3.1, we have
\begin{eqnarray}\label{s3.31}
&&\eta_h(F,\psi_{h},\varphi_{h},\Omega)\lesssim
\|(\psi,\varphi)-(\psi_{h},\varphi_{h})\|_{h}\nonumber\\
&&~~~~~~+Osc_{m}(F)+Osc_{m}(\eta_{\ell,3}(\psi_{h}))+Osc_{m}(\eta_{\ell,4}(\psi_{h})).
\end{eqnarray}

\section{The a posteriori error analysis of $C^{0}$IPG discrete scheme for the eigenvalue problem (\ref{s2.5}) }

\indent Now, we analyze the a posteriori error of the $C^{0}$IPG eigenpair $(\lambda_{h}, u_{h}, \omega_{h})$.\\
\indent Consider the source problem (\ref{s2.9}) associated with
(\ref{s2.5}) with $(f,g)=\lambda_{h}(u_{h},\omega_{h})$. Then its generalized solution
 $(\psi,\varphi)=\lambda_{h}T(u_{h},\omega_{h})$ and the $C^0$IPG approximation
$(\psi_{h},\varphi_{h})=\lambda_{h}T_{h}(u_{h},\omega_{h})=(u_{h},\omega_{h})$. Let $v=0$ in (\ref{s2.25}), we get $\omega_{h}=\lambda_{h}u_{h}$.
Thus, in (\ref{s3.6}), we have
\begin{eqnarray*}
&&(\sum\limits_{\ell\in
\mathcal{E}}\|\frac{n+1}{n-1}\|_{0,\ell}^2h_{\ell}^{4}\eta_{\ell,1}(f)^{2})^{\frac{1}{2}}=(\sum\limits_{\ell\in
\mathcal{E}}\|\frac{n+1}{n-1}\|_{0,\ell}^2h_{\ell}^{4}\eta_{\ell,1}(\lambda_{h}u_{h})^{2})^{\frac{1}{2}}\approx(\sum\limits_{\ell\in
\mathcal{E}}h_{\ell}^{4}\eta_{\ell,1}(u_{h})^2)^{\frac{1}{2}},\nonumber\\
&&h^4\sum\limits_{\ell\in
\mathcal{E}\cap\partial\kappa}\eta^2_{\ell,1}(\varphi_{h})=h^4\sum\limits_{\ell\in
\mathcal{E}\cap\partial\kappa}\eta^2_{\ell,1}(\lambda_{h}u_{h})\approx\sum\limits_{\ell\in
\mathcal{E}\cap\partial\kappa}h^4\eta^2_{\ell,1}(u_{h}),\\
&&\|f-\varphi_{h}\|_{0,\kappa}^{2}=\|\lambda_{h}u_{h}-\omega_{h}\|_{0,\kappa}^{2}=0.
\end{eqnarray*}
Hence, from (\ref{s3.6}), (\ref{s3.7}), (\ref{s3.15}) and (\ref{s3.31}) we obtain
\begin{eqnarray}
&&\eta^2_h(F,u_{h},\omega_{h},\kappa)=\eta^2_{\kappa}(F,u_{h})
+\sum_{\ell\in\mathcal{E}^{b}\cap\partial\kappa}\eta^2_{\ell,1}(u_{h})
+\frac{1}{2}\sum_{\ell\in\mathcal{E}^{i}\cap\kappa}
\{\eta^2_{\ell,1}(u_{h})\nonumber\\
&&~~~+\eta^2_{\ell,2}(u_{h})+\eta^2_{\ell,3}(u_{h})
+\eta^2_{\ell,4}(u_{h})\}
+O(\sum\limits_{\ell\in
\mathcal{E}\cap\partial\kappa}h^4\eta^2_{\ell,1}(u_{h})),\nonumber\\
&&\eta^2_h(F,u_{h},\omega_{h},\Omega)
=\sum_{\kappa\in\pi_{h}}\eta^2_h(F,u_{h},\omega_{h},\kappa),\nonumber\\\label{s4.1}
&&\|\lambda_{h}T(u_{h},\omega_{h})-\lambda_{h}T_{h}(u_{h},\omega_{h})\|_{h}\lesssim
\eta_h(F,u_{h},\omega_{h},\Omega),\\\label{s4.2}
&&\eta_h(F,u_{h},\omega_{h},\Omega)\lesssim
\|\lambda_{h}T(u_{h},\omega_{h})-\lambda_{h}T_{h}(u_{h},\omega_{h})\|_{h}\nonumber\\
&&~~~~~~+Osc_{m}(F)+Osc_{m}(\eta_{\ell,3}(u_{h}))+Osc_{m}(\eta_{\ell,4}(u_{h})).
\end{eqnarray}
where $f=\lambda_{h}u_{h}$, $g=\lambda_{h}\omega_{h}$ in $F$.\\

\indent It is noted that $O(\sum\limits_{\ell\in
\mathcal{E}\cap\partial\kappa}h^4\eta^2_{\ell,1}(u_{h}))$
is higher order small than  $\sum\limits_{\ell\in\mathcal{E}^{b}\cap\partial\kappa}\eta^2_{\ell,1}(u_{h})
+\frac{1}{2}\sum\limits_{\ell\in\mathcal{E}^{i}\cap\kappa}
\eta^2_{\ell,1}(u_{h})$, so it can be neglected in actual numerical computation.\\

\indent The following lemma is a generalization of the Lemma 9.1 in \cite{babuska}.\\
\indent{\bf Lemma 4.1.}~~ Let $(\lambda,u,\omega)$ and
$(\lambda^{*},u^{*},\omega^{*})$ be the eigenpair of (\ref{s2.5})
and (\ref{s2.13}), respectively.
 Then for any $(v,z), (v^{*},z^{*})\in
 \mathbf{H}_{h}$, when
 $B((v,z), (v^{*},z^{*}))\not=0$ it is valid that
\begin{eqnarray}\label{s4.3}
&&\frac{A_{h}((v,z), (v^{*},z^{*}))}{B((v,z),
(v^{*},z^{*}))}-\lambda =\frac{A_{h}((u,\omega)-(v,z),
(u^{*},\omega^{*})-(v^{*},z^{*}))}{B((v,z),(v^{*},z^{*}))}\nonumber\\
&&~~~~~~-\lambda \frac{B((u,\omega)-(v,z),
(u^{*},\omega^{*})-(v^{*},z^{*}))}{B((v,z),(v^{*},z^{*}))}.
\end{eqnarray}
\indent{\bf Proof.}~~ See Lemma 3.5 in \cite{yang5}.
~~~$\Box$\\

\indent Referring Lemma 4.1 in \cite{yang2} we can deduce the following theorem.\\
\indent{\bf Theorem 4.1.}~ Assume that $\lambda$ and $\lambda_h$ are
the $ith$ eigenvalues of (\ref{s2.5}) and (\ref{s2.25}),
respectively, $(u_h,\omega_h)$ is a eigenfunction corresponding to $\lambda_h$
with $\|(u_h,\omega_{h})\|_h$
$=1$, the ascent $\alpha$ of $\lambda$
is equal to 1, and assume that $R(\Omega)$ holds and $n\in W^{1,\infty}(\Omega)\cap H^{2}(\Omega)$. Let $(\bar{u}_{h},\bar{\omega}_{h})$ be the orthogonal
projection of $(u_{h},\omega_{h})$ to $ran(E_{h}^*)$ in the
sense of inner product $A_{h}(\cdot,\cdot)$, and
\begin{eqnarray}\label{s4.4}
(u_{h}^{*},\omega_{h}^{*})=\frac{(\bar{u}_{h},\bar{\omega}_{h})}{\|(\bar{u}_{h},\bar{\omega}_{h})\|_{h}}.
\end{eqnarray}
Then there exist $(u,\omega)\in ran(E)$ and
$(u^{*},\omega^{*})\in ran(E^*)$ such that $(u_{h},\omega_{h})-(u,\omega)$ and $(u_{h}^{*},\omega_{h}^{*})-(u^{*},\omega^{*})$
satisfy
(\ref{s2.40})-(\ref{s2.41}) and (\ref{s2.44})-(\ref{s2.45}) respectively,
and
\begin{eqnarray}\label{s4.5}
&&|\lambda_{h}-\lambda|\lesssim
\||(u_{h},\omega_{h})-(u,\omega)\||_{h}\||(u_{h}^{*},\omega_{h}^{*})-(u^{*},\omega^{*})\||_{h}\nonumber\\
&&~~~+\|(u_{h},\omega_{h})-(u,\omega)\|_{\mathbf{H}^{1}}\|(u_{h}^{*},\omega_{h}^{*})-(u^{*},\omega^{*})\|_{\mathbf{H}^{1}}.
\end{eqnarray}
\indent{\bf Proof.}~
 From $\alpha=1$, we know $ran(E^{*})$ is the
space of eigenfunctions associated with $\lambda^{*}$. Chose
$(u,\omega)\in ran(E)$ such that
(\ref{s2.40})-(\ref{s2.41}) hold.
Define
\begin{eqnarray*}
f((v,z))=A(E(v,z), (u,\omega)),~~~\forall (v,z)\in \mathbf{H}.
\end{eqnarray*}
Since  for all $(v,z)\in\mathbf{H}$ one has
\begin{eqnarray*}
&&|f((v,z))|=|A(E(v,z), (u,\omega))|\le
\|E(v,z)\|_A \|(u,\omega)\|_{A}\\
&&~~~ \lesssim \sqrt{\lambda}\|E(v,z)\|_{\mathbf{H}^{1}}  \lesssim
\|E\|_{\mathbf{H}^{1}}\|(v,z)\|_A,
\end{eqnarray*}
$f$ is a linear and bounded functional on $\mathbf{H}$ and
$\|f\|_{A}\lesssim \|E\|_{\mathbf{H}^{1}}$.
 Using the Riesz Theorem, we know  there exists
$(u^{*},\omega^{*})\in \mathbf{H}$ satisfying $\|(u^{*},
\omega^{*})\|_{A}=\|f\|_{A}$ and
\begin{eqnarray}\label{s4.6}
A((v,z), (u^{*}, \omega^{*}))=A(E(v,z), (u,\omega)).
\end{eqnarray}
For any $ (v,z)\in \mathbf{H}$, notice $E(I-E)(v,z)=0$, then
\begin{eqnarray*}
&&A((v,z),(\lambda^{*-1}-T^{*})(u^{*}, \omega^{*})
)=A((\lambda^{-1}-T)(v,z),(u^{*}, \omega^{*}) )\\
&&=A((\lambda^{-1}-T)E(v,z),(u^{*}, \omega^{*})
)+A((\lambda^{-1}-T)(I-E)(v,z),(u^{*}, \omega^{*}) )=0,
\end{eqnarray*}
i.e., $(\lambda^{*-1}-T^{*})(u^{*}, \omega^{*})=0$, hence $(u^{*},
\omega^{*})\in ran(E^{*})$. By (\ref{s4.6}) we have
\begin{eqnarray}\label{s4.7}
&&\lambda B((u,\omega), (u^{*}, \omega^{*}))
=A((u,\omega), (u^{*}, \omega^{*}))=A(E(u,\omega), (u,
\omega))\nonumber\\
&&~~~=A((u,\omega), (u, \omega))\approx A_{h}((u_{h},\omega_{h}), (u_{h}, \omega_{h}))\approx 1.
\end{eqnarray}
Then, there exits $(\bar{u}^*_{h},\bar{\omega}^*_{h})\in ran(E_h^*)$ such that $(\bar{u}^*_{h},\bar{\omega}^*_{h})-(u^*,\omega^*)$ satisfies (\ref{s2.42}), and from (\ref{s2.40}), (\ref{s2.42}) and (\ref{s4.7}), when $h$
is small enough, there is a positive constant $C_{0}$ independent of
$h$ such that
\begin{eqnarray*}
|B((u_{h},\omega_{h}),(\bar{u}_{h}^{*},\bar{\omega}_{h}^{*}))|\geq C_{0}.
\end{eqnarray*}
Since $(\bar{u}_{h},\bar{\omega}_{h})$ is the orthogonal projection
of $(u_{h},\omega_{h})$ to $ran(E_{h}^*)$ in the sense of inner
product $A_{h}(\cdot,\cdot)$,
\begin{eqnarray*}
&&|B((u_{h},\omega_{h}),(u_{h}^{*},\omega_{h}^{*}))|=|\frac{1}{\lambda_{h}}A_{h}((u_{h},\omega_{h}),(u_{h}^{*},\omega_{h}^{*}))|\nonumber\\
&&~~~\geq
|\frac{1}{\lambda_{h}}A_{h}((u_{h},\omega_{h}),\frac{( \bar{u}_{h}^{*},\bar{\omega}_{h}^{*})}{\|(\bar{u}_{h}^*,\bar{\omega}_{h}^*)\|_{h}})|\nonumber\\
&&~~~\geq \frac{1}{\|(\bar{u}_{h}^*,\bar{\omega}_{h}^*)\|_{h}}
|B((u_{h},\omega_{h}),(\bar{u}_{h}^{*},\bar{\omega}_{h}^{*}))|
\gtrsim C_{0}.
\end{eqnarray*}
In (\ref{s4.3}), chose $(v,z)=(u_{h},\omega_{h})$ and
$(v^{*},z^{*})=(u_{h}^{*},\omega_{h}^{*})$, and chose
$(u^{*},\omega^{*})$ such that $(u_{h}^{*},\omega_{h}^{*})-(u^{*},\omega^{*})$ satisfies
(\ref{s2.44})-(\ref{s2.45}), noting that
$$\lambda_{h}=A((u_{h},\omega_{h}),(u_{h}^{*},\omega_{h}^{*}))/B((u_{h},\omega_{h}),(u_{h}^{*},\omega_{h}^{*})),$$
we obtain (\ref{s4.5}).
~~~$\Box$\\
\indent {\bf Remark 4.1.}~~When $\lambda$ is a simple eigenvalue, $ran(E_{h}^{*})$ is a one-dimensional space
spanned by the eigenfunction $(u_{h}^{*},\omega_{h}^{*})$ of (\ref{s2.35}) with the mesh size $h$.
When the multiplicity $q>1$ of $\lambda$, in actual computation we can use the two sided Arnoldi algorithm to compute both left and right eigenfunctions of (\ref{s2.25}) at the
  same time,
  and obtain $(u_{h},\omega_{h})$ and $(u_{h}^{*},\omega_{h}^{*})$.
\\

\indent {\bf Lemma 4.2.}~~Assume that the ascent $\alpha=1$ of $\lambda$, $(u_{h},\omega_{h})$ is an eigenfunction
corresponding to $\lambda_{h}$ and $\|(u_{h},\omega_{h})\|_h=1$,
then there exists eigenfunction $(u,\omega)$ corresponding to
$\lambda$ such that
\begin{eqnarray}\label{s4.8}
|\lambda_{h}-\lambda|+\|(u_{h},\omega_{h})-(u,\omega)\|_{\mathbf{H}^1}\lesssim\|(T-T_{h})(u_{h},\omega_{h})\|_{\mathbf{H}^1}.
\end{eqnarray}
\indent {\bf Proof.}~~Using the argument as in proposition 5.3 in \cite{chatelin} we can deduce
\begin{eqnarray}\label{s4.9}
\|(u_{h},\omega_{h})-(u,\omega)\|_{\mathbf{H}^1}\lesssim\|(T-T_{h})(u_{h},\omega_{h})\|_{\mathbf{H}^1}.
\end{eqnarray}
Simple calculation shows
\begin{eqnarray*}
&&B((T-T_{h})(u_{h},\omega_{h}), (u^*,\omega^*))=B(T(u_{h},\omega_{h}), (u^*,\omega^*))-B(T_{h}(u_{h},\omega_{h}), (u^*,\omega^*))\nonumber\\
&&~~~=\lambda^{-1}A(T(u_{h},\omega_{h}), (u^*,\omega^*))-B(T_{h}(u_{h},\omega_{h}), (u^*,\omega^*))\nonumber\\
&&~~~=\lambda^{-1}B((u_{h},\omega_{h}),(u^*,\omega^*))-\lambda_{h}^{-1}B((u_{h},\omega_{h}), (u^*,\omega^*))\nonumber\\
&&~~~=(\lambda^{-1}-\lambda_{h}^{-1})B((u_{h},\omega_{h}),(u^*,\omega^*)),
\end{eqnarray*}
where $(u^*,\omega^*)$ satisfies Theorem 4.1.\\
Then the above equality implies
\begin{eqnarray}\label{s4.10}
|\lambda_{h}-\lambda|\lesssim\|(T-T_{h})(u_{h},\omega_{h})\|_{\mathbf{H}^1}.
\end{eqnarray}
Combining (\ref{s4.9}) and (\ref{s4.10}) we get (\ref{s4.8}).~~~$\Box$\\
\indent Referring \cite{dai1} et al., we give the relationship
between the $C^{0}$IPG eigenvalue approximation and the associated
$C^{0}$IPG boundary value approximation.\\
\indent {\bf Lemma 4.3.}~~
Let $(\lambda_{h},(u_{h},\omega_{h}))$ be the $ith$ eigenpair of (\ref{s2.25})
with $\|(u_{h},\omega_{h})\|_{h}=1$, $\lambda$ be the $ith$ eigenvalue of
(\ref{s2.5}), then there exists an eigenfunction $(u,\omega)$  corresponding
to $\lambda$, such that
\begin{eqnarray}\label{s4.11}
\|(u_{h},\omega_{h})-(u,\omega)\|_{h}&=&\lambda_{h}\|T(u_{h},\omega_{h})-T_{h}(u_{h},\omega_{h})\|_{h}+R_1,
\end{eqnarray}
where $\mid R_1\mid \lesssim\|(T-T_{h})(u_{h},\omega_{h})\|_{\mathbf{H}^{1}}$.\\
\indent {\bf Proof.}
From (\ref{s2.11}), (\ref{s2.12}) and (\ref{s4.8}) we have
\begin{eqnarray}\label{s4.12}
&&\|(u,\omega)-\lambda_{h}T(u_{h},\omega_{h})\|_{h}=\|\lambda T(u,\omega)-\lambda_{h}T(u_{h},\omega_{h})\|_{h}\nonumber\\
&&~~~\lesssim \|\lambda(u,\omega)-\lambda_{h}(u_{h},\omega_{h})\|_{\mathbf{H}^1}\lesssim \|(T-T_{h})(u_{h},\omega_{h})\|_{\mathbf{H}^1}.
\end{eqnarray}
Denote
\begin{eqnarray}\label{s4.13}
\|(u_{h},\omega_{h})-(u,\omega)\|_{h}= \lambda_{h}\|(T-T_{h})(u_{h},\omega_{h})\|_{h}+R_1.
\end{eqnarray}
From the triangle inequality and (\ref{s4.12}) we deduce
\begin{eqnarray}\label{s4.14}
&&\mid R_{1}\mid=\mid\|(u_{h},\omega_{h})-(u,\omega)\|_{h}- \lambda_{h}\|(T-T_{h})(u_{h},\omega_{h})\|_{h}\mid\nonumber\\
&&~~~=\mid\|(u_{h},\omega_{h})-(u,\omega)\|_{h}-\||\lambda_{h}T(u_{h},\omega_{h})-(u_{h},\omega_{h})\|_{h}\mid\nonumber\\
&&~~~\leq
\|(u,\omega)-\lambda_{h} T(u_{h},\omega_{h})\|_{h
}\lesssim\|(T-T_{h})(u_{h},\omega_{h})\|_{\mathbf{H}^{1}}.
\end{eqnarray}
Due to (\ref{s4.13}) and (\ref{s4.14}), (\ref{s4.11}) is
obtained.~~~$\Box$\\
\indent{\bf Theorem 4.2.}~~
Let $(\lambda_{h},(u_{h},\omega_{h}))$ be the $ith$ eigenpair of (\ref{s2.25})
with $\|(u_{h},\omega_{h})\|_{h}=1$, $\lambda$ be the $ith$ eigenvalue of
(\ref{s2.5}).  Assume
that $R(\Omega)$ holds and $n\in W^{1,\infty}(\Omega)\cap H^{2}(\Omega)$, then there exists an eigenfunction $(u,\omega)$  corresponding
to $\lambda$, such that
\begin{eqnarray}\label{s4.15}
&&\|(u_{h},\omega_{h})-(u,\omega)\|_{h}\lesssim\eta_h(F, u_h,\omega_{h},\Omega),\\\label{s4.16}
&&\eta_h(F, u_h,\omega_{h},\Omega)\lesssim\|(u_{h},\omega_{h})-(u,\omega)\|_{h}\nonumber\\
&&~~~~~~
+Osc_{m}(F)+Osc_{m}(\eta_{\ell,3}(u_{h}))+Osc_{m}(\eta_{\ell,4}(u_{h})).
\end{eqnarray}
\indent{\bf Proof.}
Combining  (\ref{s4.11}) with (\ref{s4.1}) we get (\ref{s4.15}).
Combining  (\ref{s4.11}) with (\ref{s4.2}) and neglecting the higher order small quantity $R_1$ we get (\ref{s4.16}). ~~~$\Box$\\

\indent For the dual problem (\ref{s2.13}),
denote
$$F^*=F^*(f,g)=\frac{-1}{n-1}\Delta f
-\Delta(\frac{n}{n-1}f)
-\frac{n}{n-1}g.$$
Using the same argument as in Theorem 4.2 we can prove the following theorem.\\
\indent{\bf Theorem 4.3.}~~
Let $(\lambda_{h}^*,(u_{h}^*,\omega^*_{h}))$ be the $ith$ eigenpair of (\ref{s2.35})
with $\|(u_{h}^*,\omega_{h}^*)\|_{h}=1$, $\lambda^*$ be the $ith$ eigenvalue of
(\ref{s2.13}).  Assume
that $R(\Omega)$ holds and $n\in W^{1,\infty}(\Omega)\cap H^{2}(\Omega)$, then there exists an eigenfunction $(u^*,\omega^*)$  corresponding
to $\lambda^*$, such that
\begin{eqnarray}\label{s4.17}
&&\|(u_{h}^*,\omega_{h}^*)-(u^*,\omega^*)\|_{h}\lesssim\eta_h(F^*, u_{h}^*,\omega_{h}^*,\Omega),\\\label{s4.18}
&&\eta_h(F^*, u_{h}^*,\omega_{h}^*,\Omega)\lesssim\|(u_{h}^*,\omega_{h}^*)-(u^*,\omega^*)\|_{h}\nonumber\\
&&~~~~~~+Osc_{m}(F^*)+Osc_{m}(\eta_{\ell,3}(u_{h}^*))+Osc_{m}(\eta_{\ell,4}(u_{h}^*)),
\end{eqnarray}
where $f=\lambda_{h}^*u_{h}^*$, $g=\lambda_{h}^*\omega_{h}^*$ in $F^*$.\\
\indent{\bf Theorem 4.4.}~~Under the condition of Theorem 4.1,
the following estimate holds
\begin{eqnarray}\label{s4.19}
&&|\lambda_{h}-\lambda|\lesssim\eta_h^2(F, u_h,\omega_{h},\Omega)+\eta_h^2(F^*, u_{h}^*,\omega_{h}^*,\Omega)+R_2,
\end{eqnarray}
where
\begin{eqnarray}
R_2=\sum\limits_{\kappa\in\pi_{h}}h_{\kappa}^{2\alpha}\|(u,\omega)-I_{h}(u,\omega)\|_{H^{2+\alpha}(\kappa)}^{2}
+\sum\limits_{\kappa\in\pi_{h}}h_{\kappa}^{2\alpha}\|(u^{*},\omega^{*})-I_{h}(u^{*},\omega^{*})\|_{H^{2+\alpha}(\kappa)}^{2}.\nonumber
\end{eqnarray}
{\bf Proof.}~~ Thanks to Poincar$\acute{e}$-Friedrichs
inequalities in \cite{brenner4}, we have
$\|(u_{h},\omega_{h})-(u,\omega)\|_{\mathbf{H}^{1}}\lesssim\||(u_{h},\omega_{h})-(u,\omega)\||_{h}$
and
$\|(u_{h}^{*},\omega_{h}^{*})-(u^{*},\omega^{*})\|_{\mathbf{H}^{1}}\lesssim\||(u_{h}^{*},\omega_{h}^{*})-(u^{*},\omega^{*})\||_{h}$.
Thus from (\ref{s4.5}) we get
\begin{eqnarray}\label{s4.20}
&&|\lambda_{h}-\lambda|\lesssim
\||(u_{h},\omega_{h})-(u,\omega)\||_{h}\||(u_{h}^{*},\omega_{h}^{*})-(u^{*},\omega^{*})\||_{h}
\end{eqnarray}
Due to (\ref{s2.29}), the triangle inequality, (\ref{s2.30}), (\ref{s2.31}) and the interpolation estimate, we deduce
\begin{eqnarray*}
&&\||(u_{h},\omega_{h})-(u,\omega)\||^2_h\nonumber\\
&&~~~\leq
(\||(u_{h},\omega_{h})-I_{h}(u,\omega)\||+\||(u,\omega)-I_{h}(u,\omega)\||_h)^{2}\nonumber\\
&&~~~\lesssim (\|(u_{h},\omega_{h})-I_{h}(u,\omega)\|_h+\||(u,\omega)-I_{h}(u,\omega)\||_h)^{2}\nonumber\\
&&~~~\lesssim
(\|(u_{h},\omega_{h})-(u,\omega)\|_h+\||(u,\omega)-I_{h}(u,\omega)\||_h)^{2}\nonumber\\
&&~~~\lesssim \eta_h^2(F, u_h,\omega_{h},\Omega)+\sum\limits_{\kappa\in\pi_{h}}h_{\kappa}^{2\alpha}\|(u,\omega)-I_{h}(u,\omega)\|_{H^{2+\alpha}(\kappa)}^{2}.\nonumber
\end{eqnarray*}
Similarly, we can get
\begin{eqnarray*}
&&\||(u_{h}^{*},\omega_{h}^{*})-(u^{*},\omega^{*})\||^2_h\nonumber\\
&&~~~\lesssim \eta_h^2(F, u_h^{*},\omega_{h}^{*},\Omega^{*})+\sum\limits_{\kappa\in\pi_{h}}h_{\kappa}^{2\alpha}\|(u^{*},\omega^{*})-I_{h}(u^{*},\omega^{*})\|_{H^{2+\alpha}(\kappa)}^{2}.\nonumber
\end{eqnarray*}
Submitting the above two estimate into (\ref{s4.20}), we get (\ref{s4.19}).~~~$\Box$\\
{\bf Remark 4.2.}~~
From Theorems 4.2 and 4.3, we know the indicator $\eta_h^2(F, u_h,\omega_h,\Omega)+\eta_h^2(F^{*}, u_h^{*},\omega_{h}^{*},\Omega)$ of the eigenfunction error $\|(u_{h},\omega_{h})-(u,\omega)\|_h^2+ \|(u_{h}^{*},\omega_{h}^{*})-(u^{*},\omega^{*})\|_h^2$ is reliable and efficient up to data oscillation, so Algorithm 1 can generate a good graded mesh, which makes approximation eigenfunctions can get the optimal convergent rate $h^{m-1}$ in $\|\cdot\|_h$.
And thus we are able to expect to  get $R_2\lesssim h^{2(m-1)}$,
thereby from (\ref{s4.19}) have
 $|\lambda_h-\lambda|\lesssim h^{2(m-1)}$.
Therefore, we think that $\eta_h^2(F, u_h,\omega_h,\Omega)+\eta_h^2(F^{*}, u_h^{*},\omega_{h}^{*},\Omega)$ can be viewed as the indicator of $\lambda_h$.
The numerical experiments in Section 5 show this indicator of $\lambda_h$ is reliable and efficient. And $\lambda_h$ can achieve the optimal convergent rate.

\section{Adaptive algorithms and Numerical Experiment}
Using the a
posteriori error estimates and consulting the existing standard algorithms (see, e.g., \cite{han,dai1} ), we present the following algorithm.\\
\indent{\bf Algorithm 1}\\
\indent Choose the parameter $\sigma, \mu, 0<\theta<1$.\\
\indent{\bf Step 1.} Set $l=0$ and pick any initial mesh $\pi_{h_l}$ with the mesh size $h_{l}$.\\
\indent{\bf Step 2.} Solve (\ref{s2.25}) on $\pi_{h_l}$ for discrete solution $(\lambda_{h_l},(u_{h_l},\omega_{h_l}))$ with
$\|(u_{h_l},\omega_{h_l})\|_{h}$
$=1$ and find $(u^*_{h_l},\omega^*_{h_l})\in ran(E_{h}^*)$
 by (\ref{s4.4}) (also see Remark 4.1).\\
\indent{\bf Step 3.} Compute the local indicators $\eta_{h_l}^2(F, u_{h_l},\omega_{h_l},\kappa)+\eta_{h_l}^2(F^{*}, u_{h_l}^{*},\omega_{h_l}^{*},\kappa)$.\\
\indent{\bf Step 4.} Construct $\hat{\pi}_{h_l}\in \pi_{h_l}$ by {\bf Marking strategy E}.\\
\indent{\bf Step 5.} Refine $\pi_{h_l}$ to get a new mesh $\pi_{h_{l+1}}$ by procedure {\bf Refine}.\\
\indent{\bf Step 6.} Set $l=l+1$ and goto Step 2.\\
\indent{\bf Marking Strategy E}\\
\indent Given parameter $0<\theta<1$:\\
\indent{\bf Step 1.} Construct a minimal subset
$\widehat{\pi}_{h_{l}}$ of $\pi_{h_{l}}$ by selecting some elements
in $\pi_{h_{l}}$ such that
\begin{eqnarray*}
&&\sum\limits_{\kappa\in
\widehat{\pi}_{h_{l}}}(\eta_{h_l}^2(F, u_{h_l},\omega_{h_l},\kappa)+\eta_{h_l}^2(F^{*}, u_{h_l}^{*},\omega_{h_l}^{*},\kappa))\nonumber\\
&&~~~\geq \theta(\eta_{h_l}^2(F, u_{h_l},\omega_{h_l},\Omega)+\eta_{h_l}^2(F^{*}, u_{h_l}^{*},\omega_{h_l}^{*},\Omega)).
\end{eqnarray*}
\indent{\bf Step 2.} Mark all the elements in
$\widehat{\pi}_{h_{l}}$.\\

We compute the transmission eigenvalues on  the unit square domain with a slit $[0,1]^2\setminus[0.5,1]$ and the L-shaped domain $[-1,1]^2\setminus[0,1]\times[-1,0]$ using Algorithm 1 with $m=2,3$. All the initial meshes are made up of congruent triangles. And the mesh sizes take $h_0=\frac{\sqrt{2}}{32}$ and $h_0=\frac{\sqrt{2}}{16}$ for the domain with a slit and the L-shaped domain, respectively. $\theta=0.25$ and $\theta=0.5$ for $m=2$ and $m=3$, respectively. Our programs uses MATLAB2012a and the iFEM package (see \cite{chen}) on a HP-Z230 workstation(CPU 3.6GHZ and RAM 32GB).\\
\indent We use the sparse solver $eigs$ to solve
(\ref{s2.25}) and (\ref{s2.35}) for eigenvalues. Before showing the
results, some symbols need to be explained:\\
$k_{j}=\sqrt{\lambda_{j}}$;\\
$\lambda_{j,h_{l}}$: the $jth$ eigenvalue derived from the $l$th iteration using Algorithm 1,
$k_{j,h_{l}}=\sqrt{\lambda_{j,h_{l}}}$;\\
$DOF$: the number of degrees of freedom.\\
The accurate eigenvalues for the problems on the two above domains are unknown. For the domain with a slit, we take $k_1\approx2.80677803, k_2\approx 2.98066000$ for $n=16$, and take $k_1\approx4.14438323, k_7\approx5.57000885-1.31142340i$ for $n=8+x-y$. For the L-shaped domain, we take $k_1\approx1.47609911, k_2\approx1.56972499$ for $n=16$, and take $k_1\approx2.30212024, k_5\approx2.92423162-0.56458999i$ for $n=8+x-y$.
All of them are obtained by Algorithm 1. And we think them relatively accurate. By computation we also know that the first ten smallest eigenvalues are all simple.\\
We present some adaptive refined mesh in Figure 1, and the curves of the error of the numerical eigenvalues in Figures $2\sim5$.\\
From Figure 1, we can see that the singularities of the eigenfunctions for the two domain mainly center on the corner points.\\
From Figures $2\sim5$, we see that the curves of the indicator are parallel to the curves of the error of $\lambda_{j,h}$, which shows the posteriori error estimators are reliable and efficient for all the cases; we also see that the accuracy of the numerical eigenvalues on adaptive meshes, better than that on uniform meshes, can get the optimal convergence order  $O(DOF^{-m+1}),~m=2,3$.\\
However, from Figures $2\sim5$, we also see that there exists fluctuation in the results on adaptive meshes when $DOF$ is large enough. This is probably the consequence of the performance of linear algebra routine on this problem. To treat such problems to get higher accurate approximation much more careful design of the routine is needed.
\begin{figure*}[htbp]
\centering
\includegraphics[width=6.8cm]{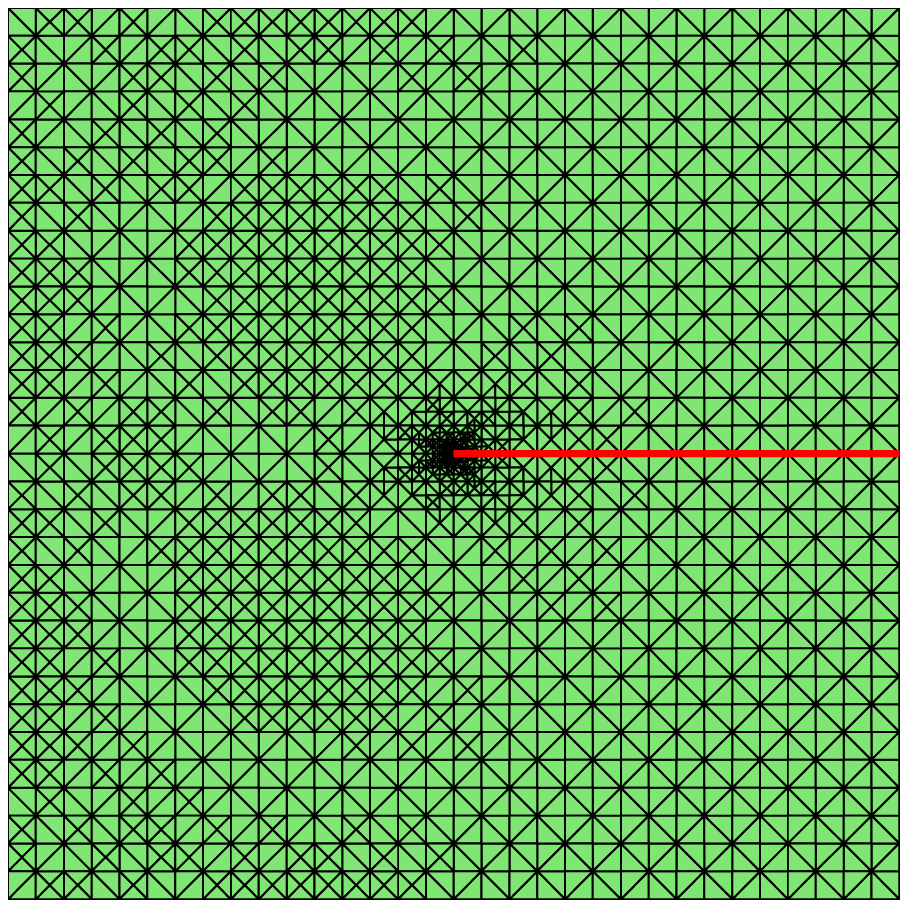}
\includegraphics[width=6.8cm]{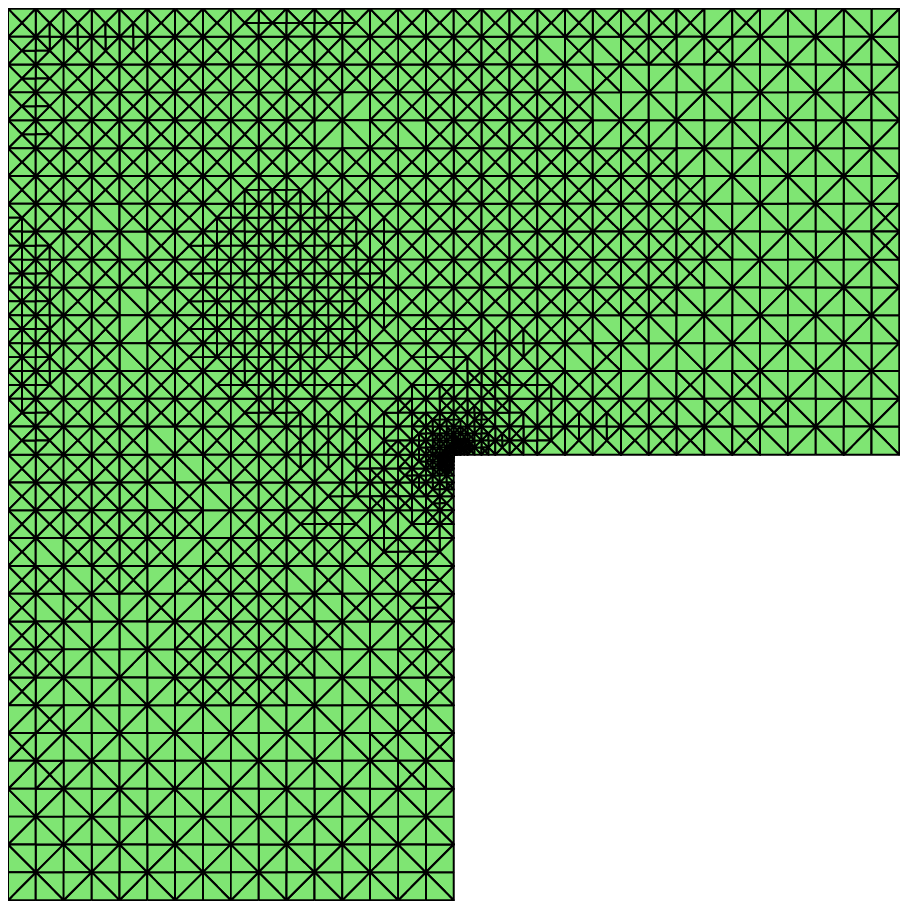};
\includegraphics[width=6.8cm]{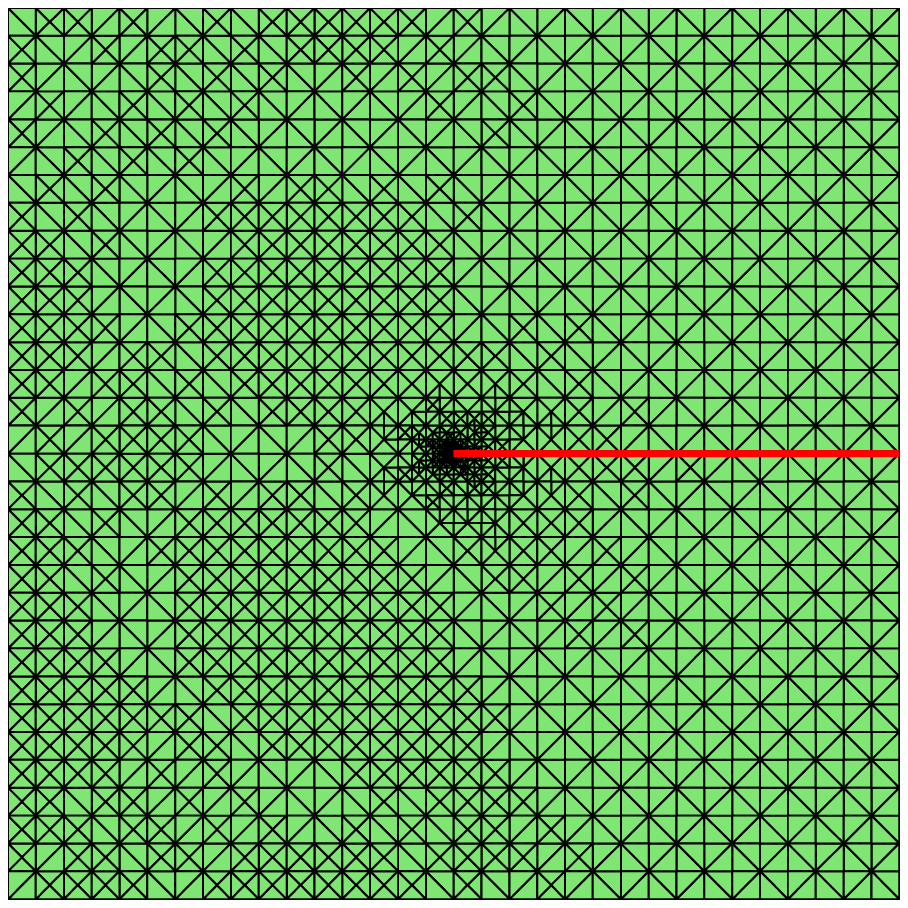}
\includegraphics[width=6.8cm]{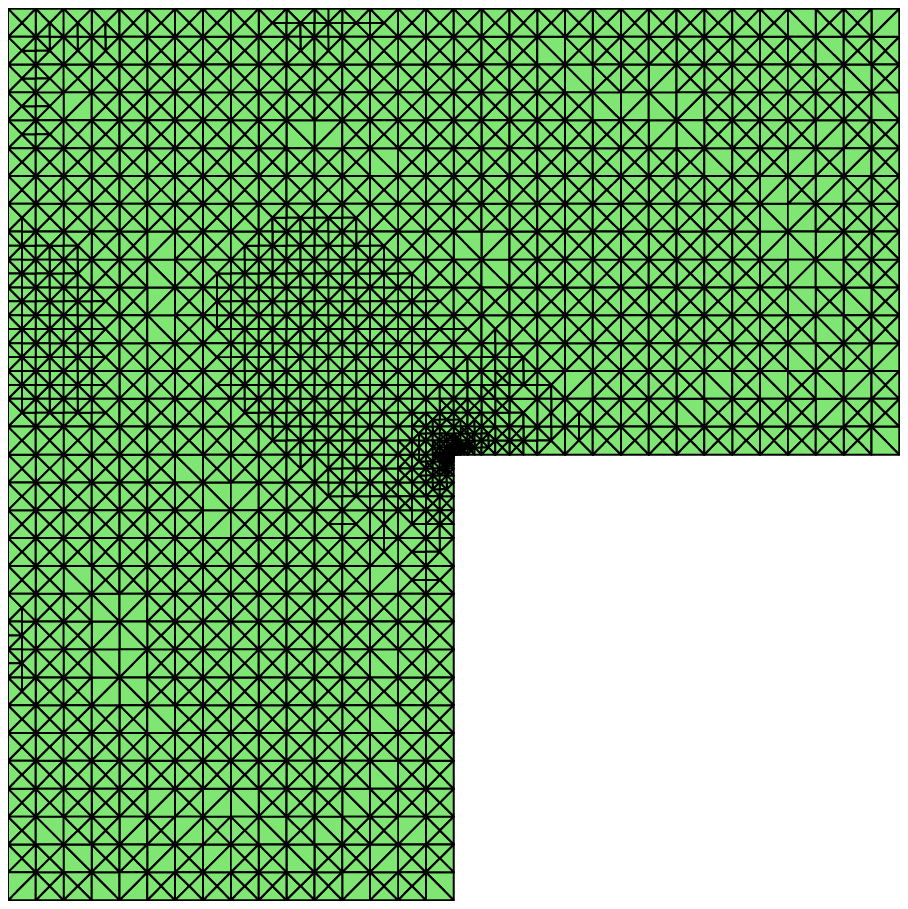}
\caption{Adaptive meshes for the smallest eigenvalue on the
domain with a slit
with $DOF=28688,n = 16, m=3$ (left up),
on the L-shaped domain with
$DOF=29972,n = 16, m=3$ (right up),
on the
domain with a slit with
$DOF=29564,n = 8+x-y, m=3$ (left down) and
on the L-shaped domain with
$DOF=32954,n = 8+x-y, m=3$ (right down).}
\end{figure*}

\begin{figure*}[htbp]
\centering
\includegraphics[width=6.8cm]{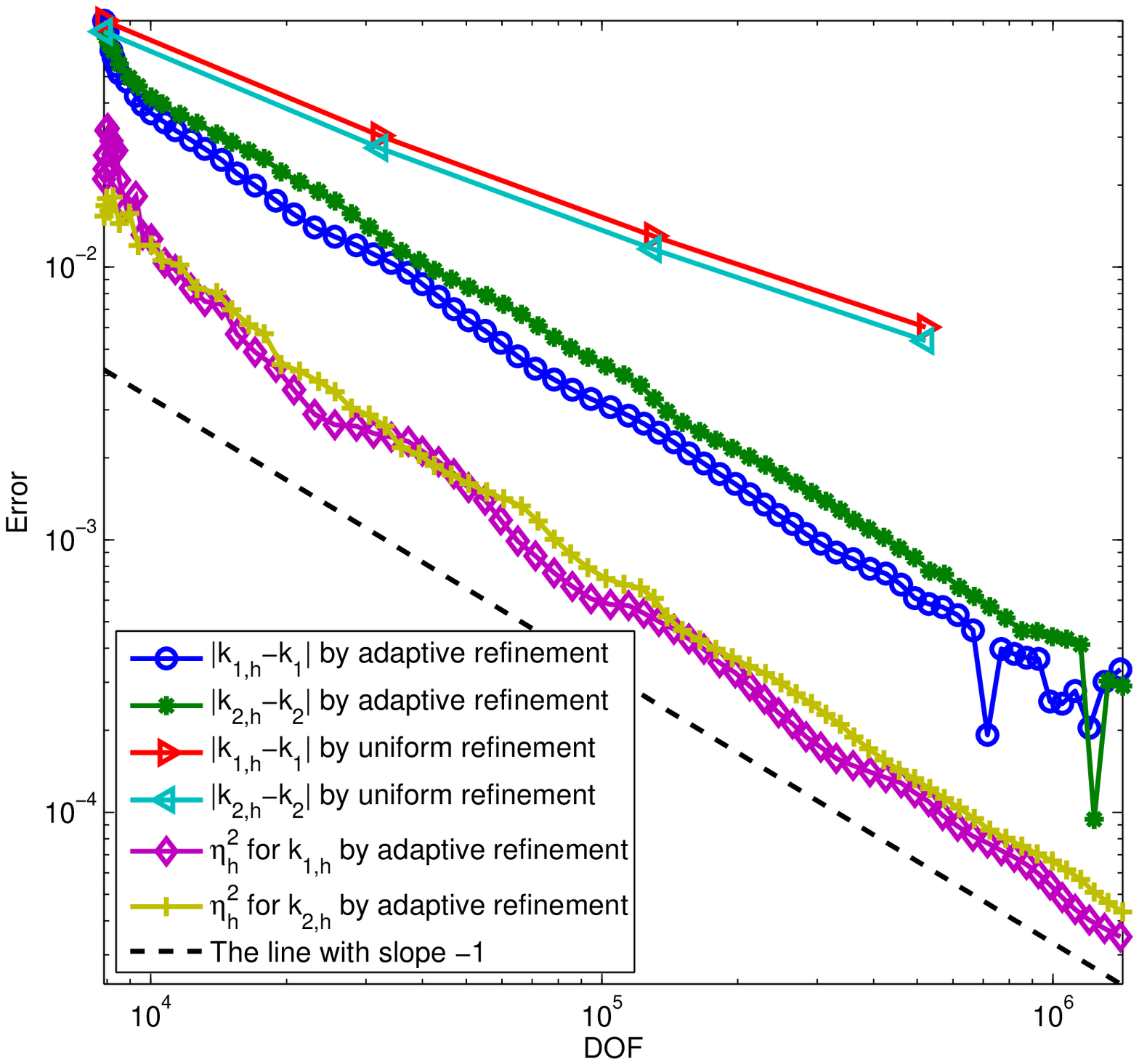}
\includegraphics[width=6.8cm]{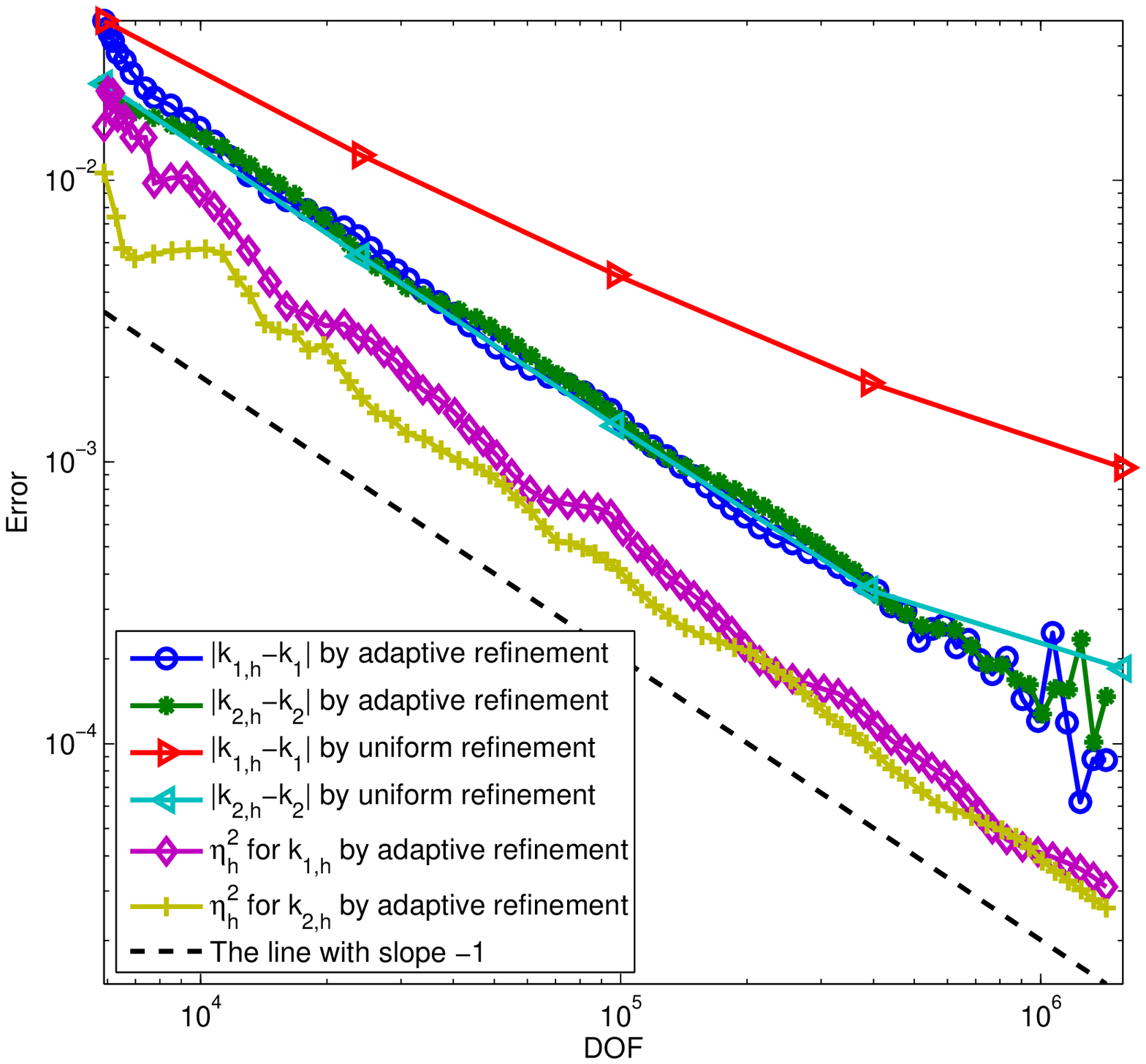}
\caption{The convergence rates of eigenvalues for the domain with a slit (left) and  for the L-shaped domain (right) when $n=16,m=2,\sigma=30,\mu=\frac{1}{15}$.}
\end{figure*}

\begin{figure*}[htbp]
\centering
\includegraphics[width=6.7cm]{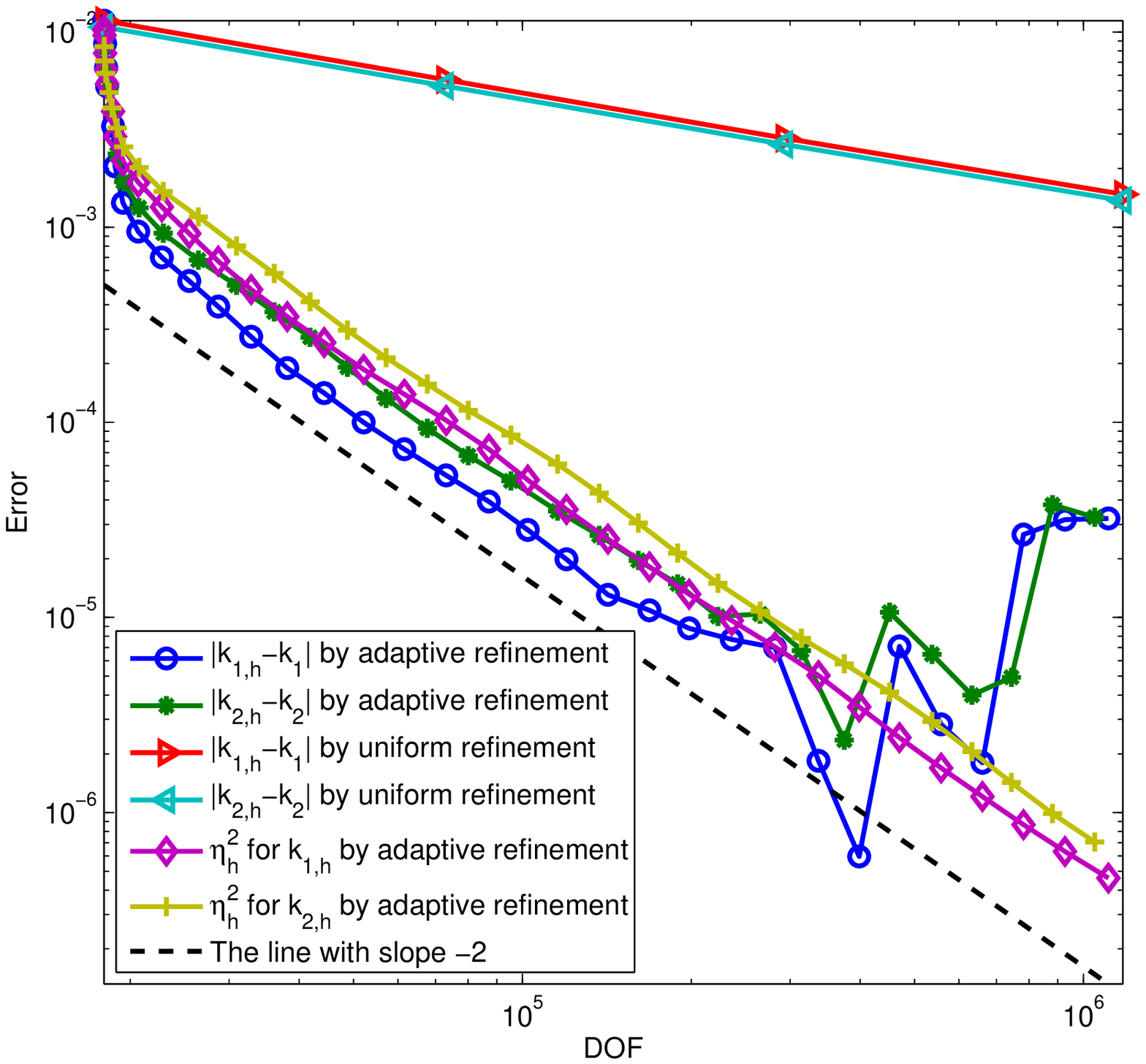}
\includegraphics[width=6.7cm]{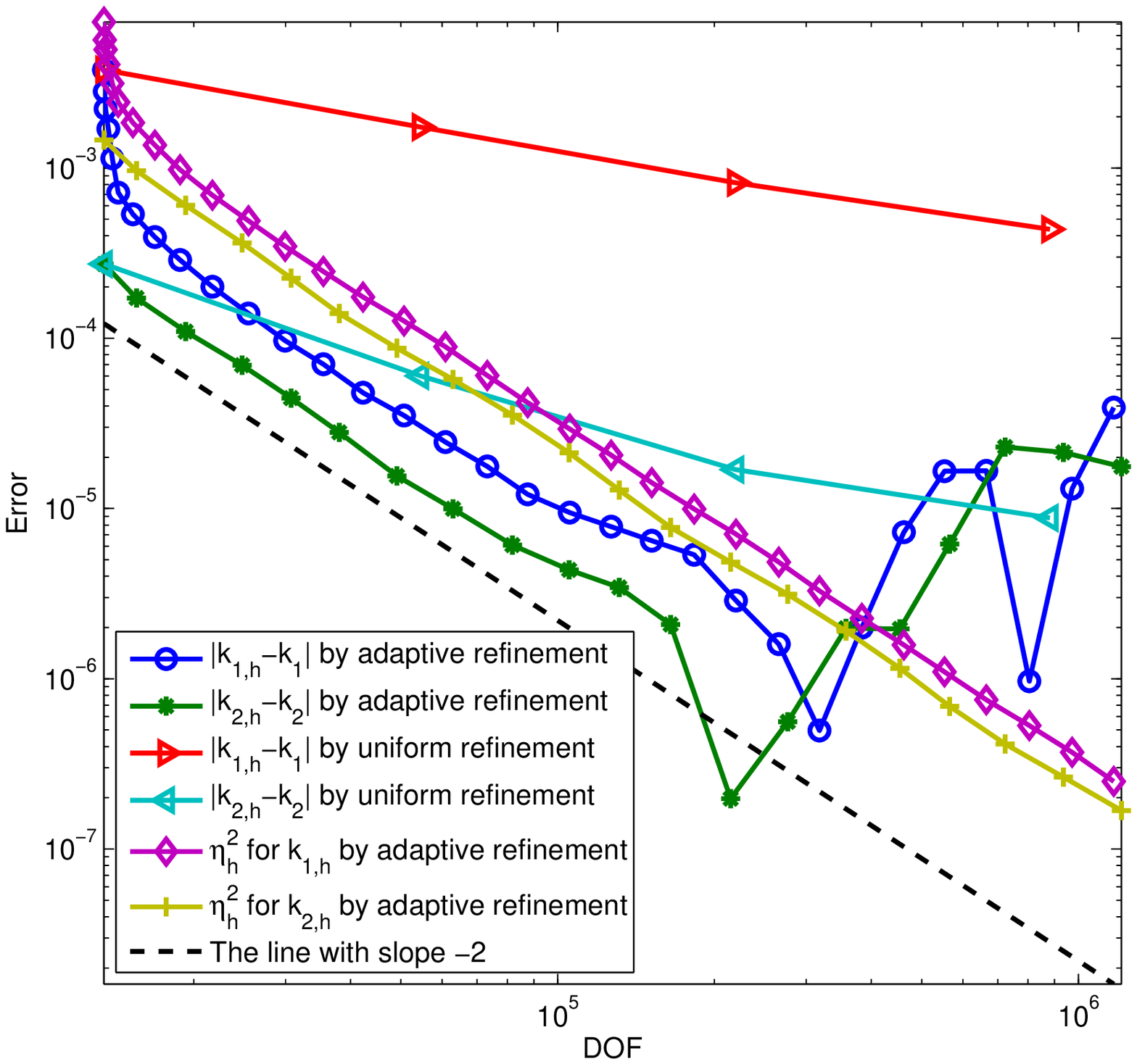}
\caption{The convergence rates of eigenvalues for the domain with a slit (left) and  for the L-shaped domain (right) when $n=16,m=3,\sigma=30,\mu=\frac{1}{15}$.}
\end{figure*}

\begin{figure*}[htbp]
\centering
\includegraphics[width=6.7cm]{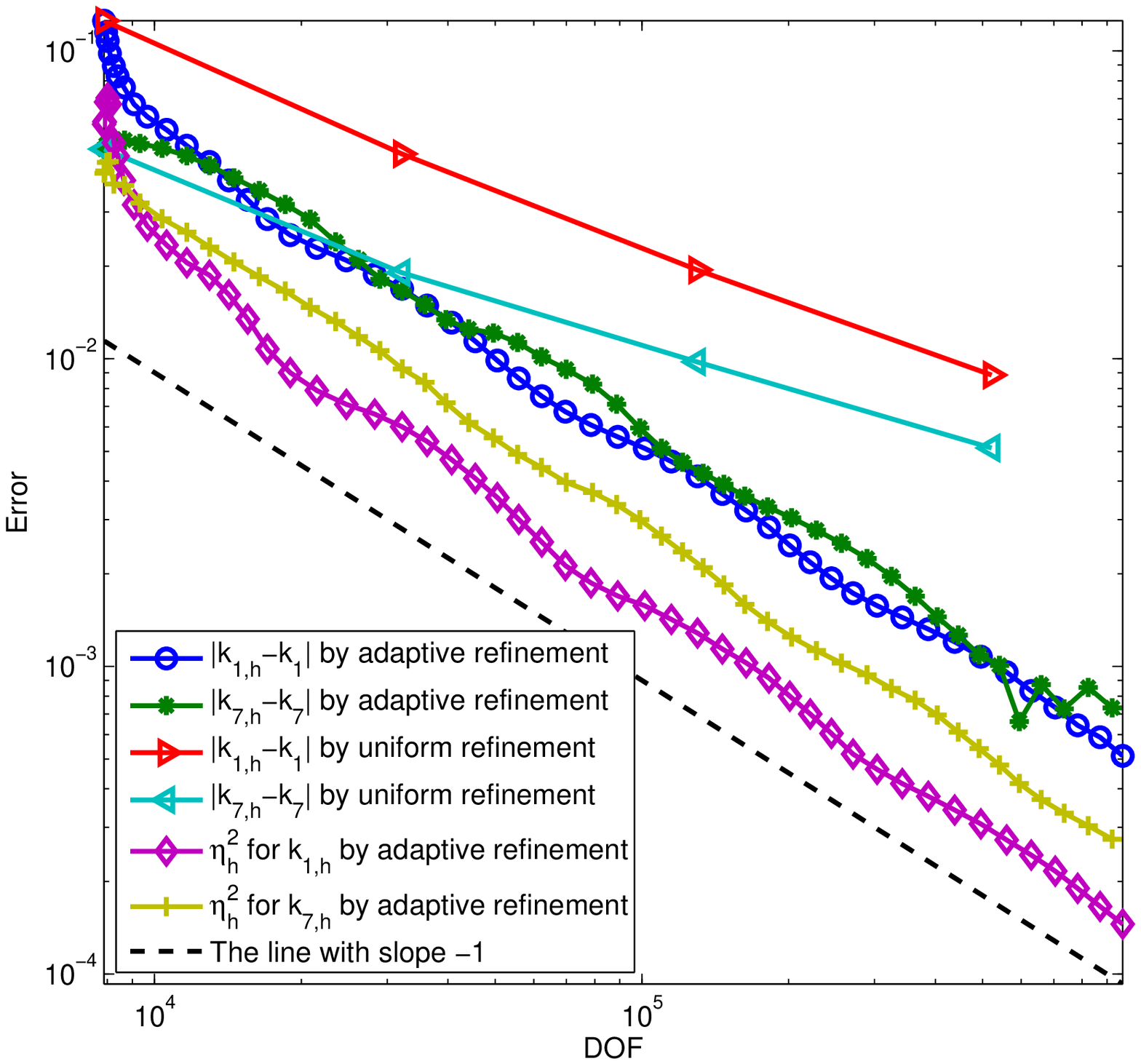}
\includegraphics[width=6.7cm]{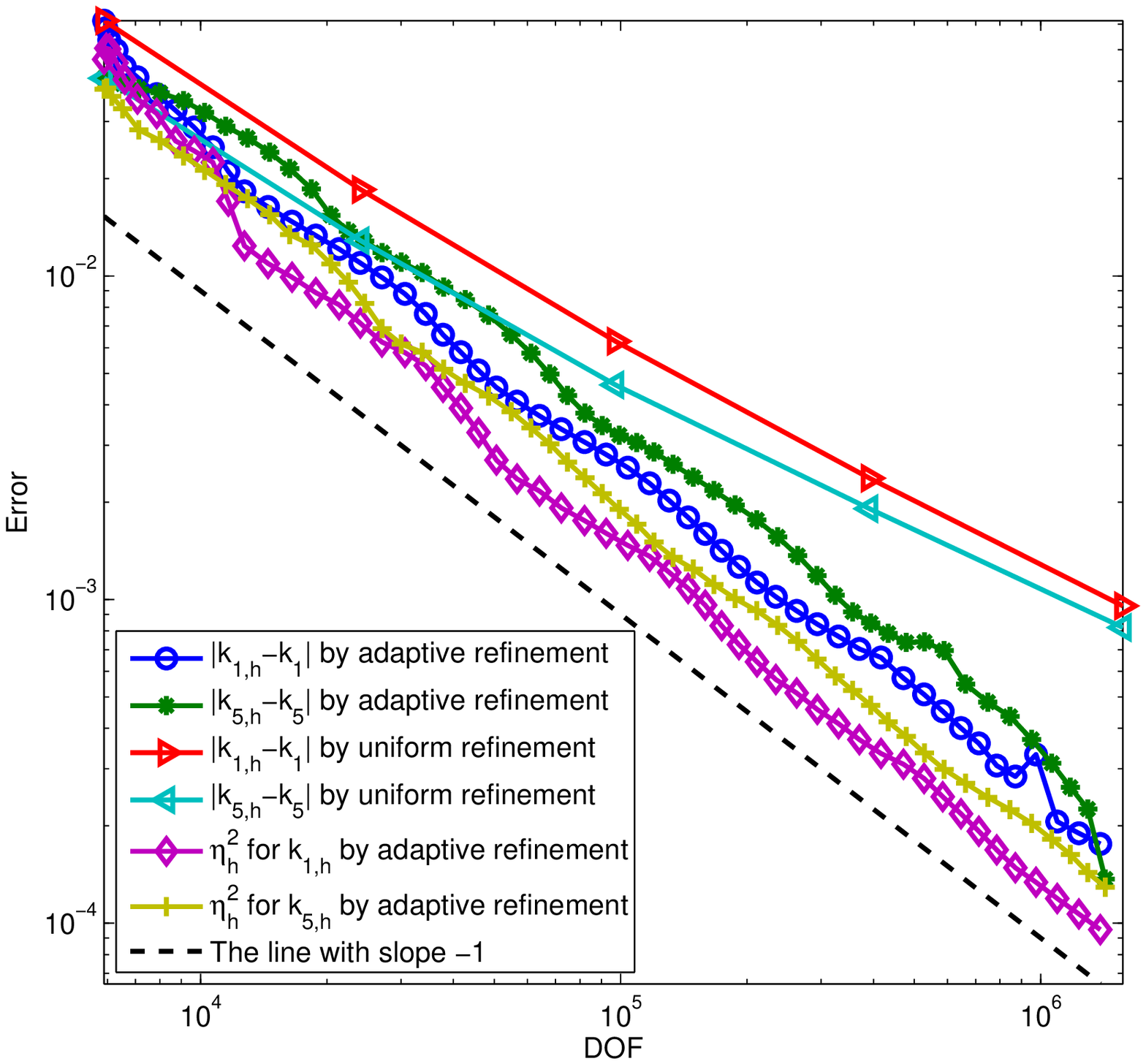}
\caption{The convergence rates of eigenvalues for the domain with a slit (left) and  for the L-shaped domain (right) when $n=8+x-y,m=2,\sigma=20,\mu=\frac{1}{9}$.}
\end{figure*}

\begin{figure*}[htbp]
\centering
\includegraphics[width=6.7cm]{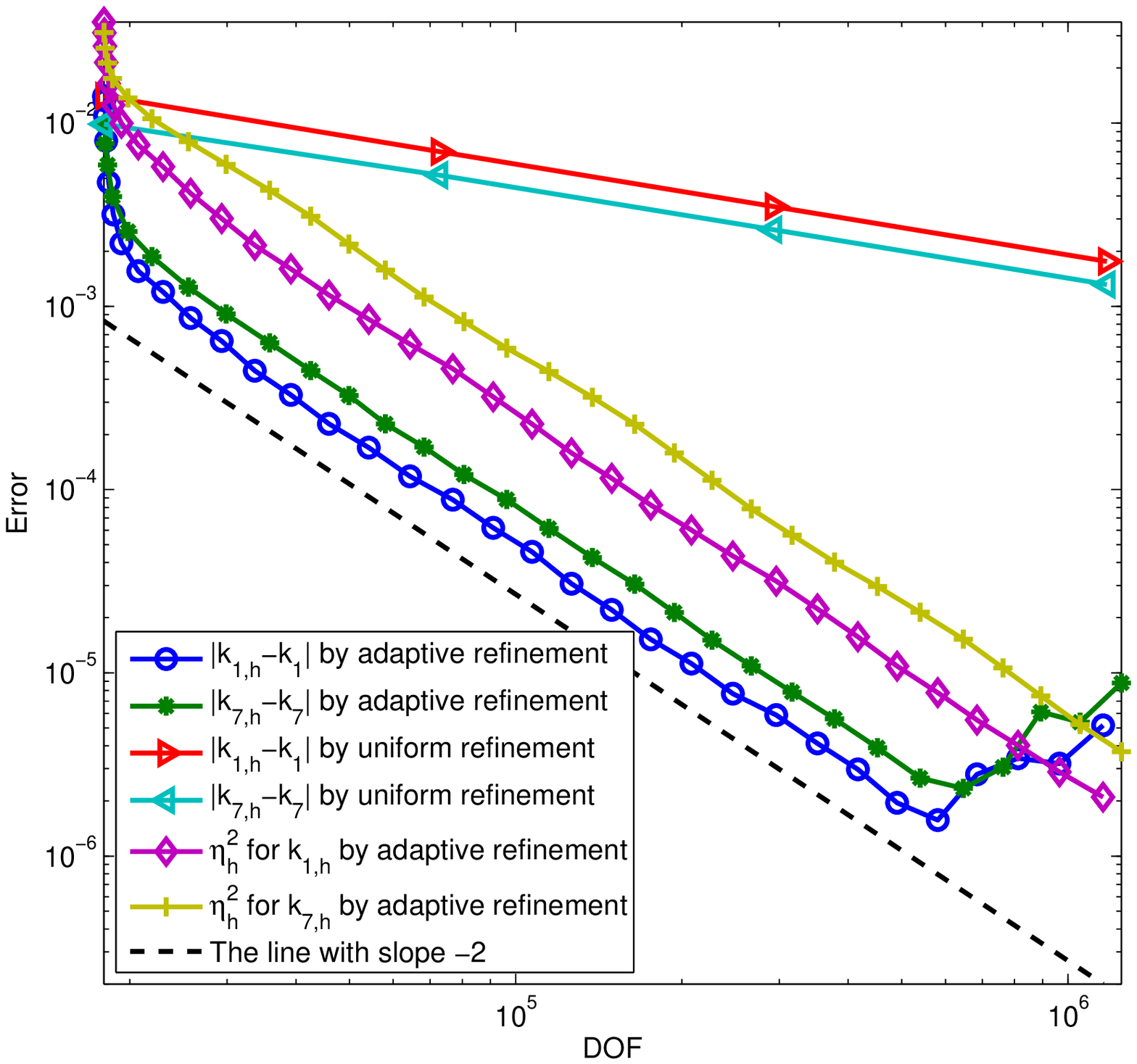}
\includegraphics[width=6.7cm]{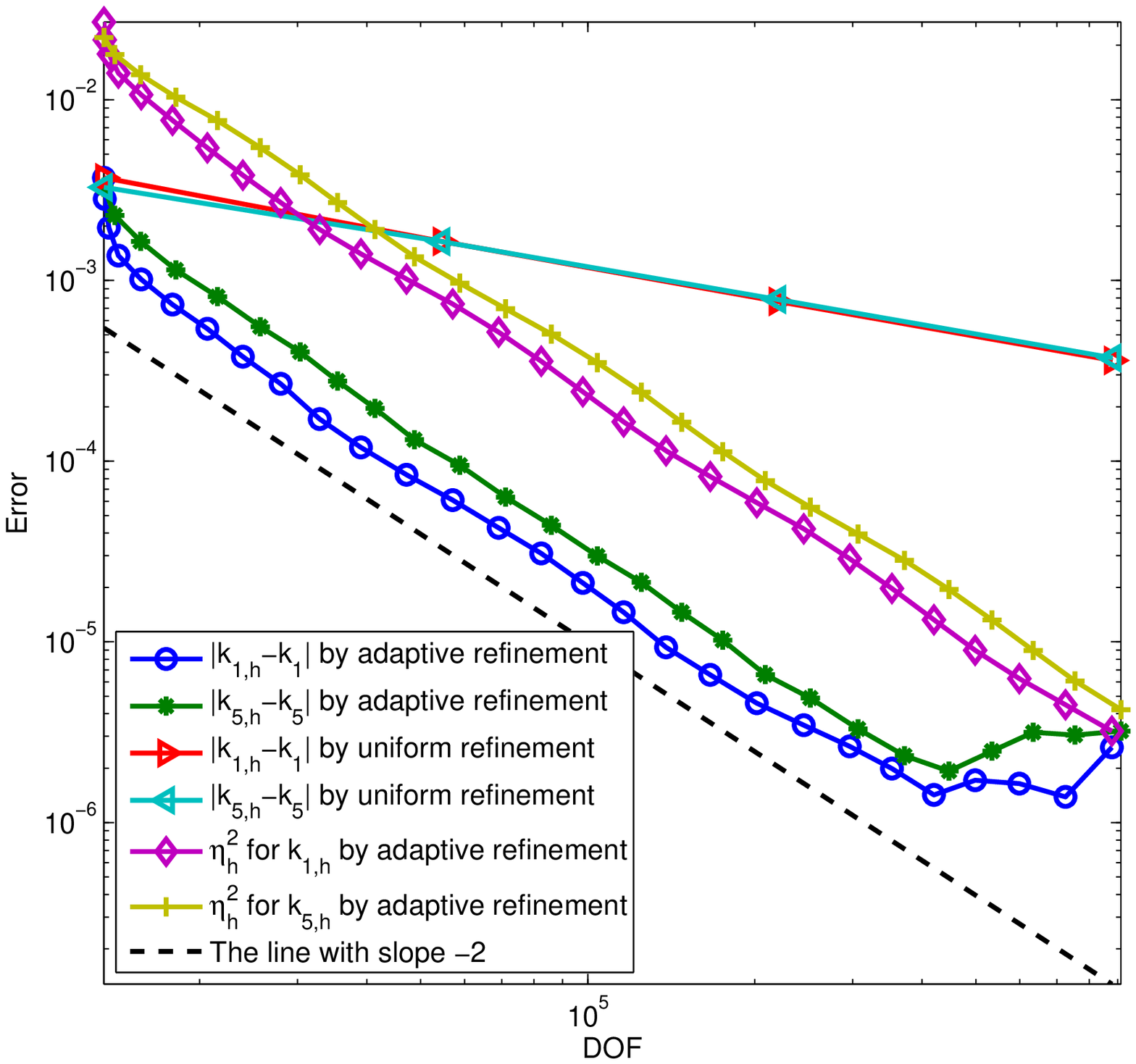}
\caption{The convergence rates of eigenvalues for the domain with a slit (left) and  for the L-shaped domain (right) when $n=8+x-y,m=3,\sigma=20,\mu=\frac{1}{9}$.}
\end{figure*}


\end{document}